\documentclass[11pt,a4paper,openany,oneside]{article}
\usepackage[a4paper,left=2.5cm,right=2.5cm, top=3cm, bottom=3cm]{geometry}
\usepackage{times}
\usepackage{authblk}
\usepackage[latin1]{inputenc}
\usepackage[english]{babel}
\usepackage{upgreek} 
\usepackage{mathrsfs}
\usepackage{stmaryrd}
\usepackage{amssymb,amsmath,amsthm}
\usepackage{hyperref}
\usepackage{graphicx}
\usepackage{epstopdf}
\usepackage{subcaption}
\usepackage{eurosym}
\usepackage{color}
\usepackage{enumerate}
\usepackage{multirow}
\usepackage{stmaryrd}
\usepackage{mathtools}
\usepackage{scalerel}
\usepackage{booktabs}
\usepackage{float}
\usepackage{amsfonts}
\usepackage{cases}
\usepackage{caption}

\theoremstyle{plain}
\newtheorem{thm}{Theorem}[section]

\newtheorem{lem}[thm]{Lemma}

\newtheorem{prop}[thm]{Proposition}

\theoremstyle{definition}
\newtheorem{defn}{Definition}[section]
\theoremstyle{remark}
\newtheorem{remark}{Remark}

\newcommand{\C}{\mathbb C}
\newcommand{\R}{\mathbb R}
\newcommand{\N}{\mathbb N}

\newcommand{\EA}{\cE^{(0)}_{ans}}
\newcommand{\Lp}{\Lambda_p}
\newcommand{\lp}{\lambda_p}
\newcommand{\La}{\Lambda_\alpha}
\newcommand{\tLp}{\widetilde\Lambda_p}
\newcommand{\tlp}{\widetilde\lambda_p}
\newcommand{\tLa}{\widetilde\Lambda_\alpha}
\newcommand{\les}{\lesssim}
\newcommand{\phk}{\partial_k^h}
\newcommand{\Zth}{Z_{th}}
\newcommand{\intpiu}{\int_{\R^n_+}}
\newcommand{\Aou}{U_{0,1}}
\newcommand{\Aod}{U_{0,2}}

\def\e{{\rm e}}
\def\ri{{\rm i}}

\newcommand{\cA}{{\mathcal A}}

\newcommand{\cD}{{\mathcal D}}
\newcommand{\cE}{{\mathcal E}}

\newcommand{\cH}{{\mathcal H}}

\newcommand{\cJ}{{\mathcal J}}

\newcommand{\cL}{{\mathcal L}}

\newcommand{\cO}{{\mathcal O}}
\newcommand{\cP}{{\mathcal P}}

\newcommand{\cT}{{\mathcal T}}
\newcommand{\cU}{{\mathcal U}}
\newcommand{\cV}{{\mathcal V}}

\DeclareFontFamily{U}{mathb}{\hyphenchar\font45}
\DeclareFontShape{U}{mathb}{m}{n}{ <-6> matha5 <6-7> matha6 <7-8>
	mathb7 <8-9> mathb8 <9-10> mathb9 <10-12> mathb10 <12-> mathb12 }{}
\DeclareSymbolFont{mathb}{U}{mathb}{m}{n}

\DeclareMathAccent{\abxring}{0}{mathb}{"38}

\DeclareFontFamily{U}{mathb}{\hyphenchar\font45}
\DeclareFontShape{U}{mathb}{m}{n}{ <-6> matha5 <6-7> matha6 <7-8>
	mathb7 <8-9> mathb8 <9-10> mathb9 <10-12> mathb10 <12-> mathb12 }{}
\DeclareSymbolFont{mathb}{U}{mathb}{m}{n}

\newcommand{\fD}{\mathfrak D}
\newcommand{\Do}{\abxring{D}}

\newcommand{\trans}{{\mathsf T}}

\DeclareMathOperator{\sign}{sign}

\DeclareMathOperator\Hdiv{H_{div}}

\DeclareOldFontCommand{\it}{\normalfont\itshape}{\mathit}

\newcommand{\dist}{\text{\rm dist}}
\newcommand{\supp}{\text{\rm supp}}

\newcommand{\iii}[1]{{\left\vert\kern-0.25ex\left\vert\kern-0.25ex\left\vert #1 \right\vert\kern-0.25ex\right\vert\kern-0.25ex\right\vert}}

\newcommand{\ccdot}{\!\cdot\!}

\numberwithin{equation}{section}

\begin{document}

\title{A quasilinear transmission problem with application to Maxwell equations with a divergence-free $\cD$-field}
\date{\today}
\author{Tom\'a\v{s} Dohnal, Giulio Romani, and Daniel P. Tietz}
\affil{Institut f\"{u}r Mathematik,  Martin-Luther-Universit\"{a}t Halle-Wittenberg,\\ 06099 Halle (Saale), Germany\\
\small{tomas.dohnal@mathematik.uni-halle.de, giulio.romani@mathematik.uni-halle.de, daniel.tietz@mathematik.uni-halle.de}}

\maketitle


\abstract{Maxwell equations in the absence of free charges require initial data with a divergence free displacement field $\cD$. In materials in which the dependence $\cD=\cD(\cE)$ is nonlinear the quasilinear problem $\nabla\cdot \cD(\cE)=0$ is hence to be solved. In many applications, e.g. in the modelling of wave-packets, an approximative asymptotic ansatz of the electric field $\cE$ is used, which satisfies this divergence condition at $t=0$ only up to a small residual. We search then for a small correction of the ansatz to enforce $\nabla\cdot \cD(\cE)=0$ at $t=0$ and choose this correction in the form of a gradient field. In the usual case of a power type nonlinearity in $\cD(\cE)$ this leads to the sum of the Laplace and $p$-Laplace operators. We also allow for the medium to consist of two different materials so that a transmission problem across an interface is produced. We prove the existence of the correction term for a general class of nonlinearities and provide regularity estimates for its derivatives, independent of the $L^2$-norm of the original ansatz. In this way, when applied to the wave-packet setting, the correction term is indeed asymptotically smaller than the original ansatz. We also provide numerical experiments to support our analysis.}

\section{Introduction}
This paper is concerned with the study of the existence and regularity estimates for weak solutions of the second-order quasilinear problem defined on the whole space $\R^n$, $n\geq 2$, 
\begin{equation}\label{eq_quasi_gen_D}
\nabla\cdot\cD(x,\nabla\phi+U_0(x))=0,\qquad x\in\R^n,
\end{equation}
where $\nabla\cdot$ stands for the divergence operator, $U_0:\R^n\to\R^n$ is a prescribed vector field and the function $\cD:\R^n\times\R^n\to\R^n$ is such that $\cD(x,\cdot)$ is growing both linearly and nonlinearly at infinity, namely
\begin{equation}\label{relationD(E)}
\cD(x,\cE):=\epsilon_1(x)\cE+\epsilon_f(x)f(\cE)\qquad\mbox{for}\;\;\cE:\R^n\to\R^n,
\end{equation}
with positive coefficients $\epsilon_1,\epsilon_f:\R^n\to\R^+$ and where $f:\R^n\to\R^n$ is a superlinear vector field, the precise behaviour of which will be specified in Section \ref{Section_Var_Ass_Funct}. However, the model growth that one should have in mind is $f(v)=|v|^{p-1}v$ for some $p>1$. Our analysis is confined to the case of positive $\epsilon_1$, $\epsilon_f$. An essential feature of our problem is that the coefficients $\epsilon_1,\epsilon_f$, although smooth in each half-spaces $\R^n_\pm:=\{x\in\R^n\,|\,\pm x_1>0\}$, are not continuous on $\Gamma:=\{x \in \R^n \,|\, x_1=0\}$. This means that \eqref{eq_quasi_gen_D} has to be understood pointwise in $\R^n_\pm$ and becomes a transmission problem with the straight interface $\Gamma$. which has to be coupled with the interface condition 
\begin{equation}\label{IFCs}
\llbracket \cD(\cdot,\nabla\phi+U_0)\cdot e_1\rrbracket(x)=0,\qquad x\in\Gamma.
\end{equation}
Here, the symbol $\llbracket\varphi\rrbracket$ denotes the jump of $\varphi$ on $\Gamma$. More precisely, for $x\in\Gamma$ we define
$$\llbracket\varphi\rrbracket (x):=\lim_{h\to0^+}\left(\varphi(x+he_1)-\varphi(x-he_1)\right),$$
where $e_1:=(1,0,\dots,0)^\trans$ and similarly we define $e_k$ for $k\in\{2,\dots,n\}$. We usually omit the variable $x$ in the notation $\llbracket\varphi\rrbracket(x)$. In Sobolev spaces the interface condition has to be understood in a trace sense, see Definition \ref{def_trace}.

\paragraph{Motivation: a Maxwell problem} Problem \eqref{eq_quasi_gen_D}-\eqref{relationD(E)}-\eqref{IFCs} is motivated by the study of Maxwell equations when one investigates the configuration of two materials separated by a straight interface. Let us illustrate this background in more detail. Maxwell equations in nonlinear dielectric materials without free charges and currents are described by
\begin{equation}\label{Max}
\upmu_0\partial_t\cH=-\nabla\times\cE,\qquad\upepsilon_0\partial_t\cD=\nabla\times\cH,\qquad\nabla\cdot\cD=\nabla\cdot\cH=0, \quad x\in\R^n, t>0
\end{equation}
with\footnote{In the two-dimensional case $n=2$ we assume without loss of generality $(\cE,\cH,\cD)=(\cE,\cH,\cD)(x_1,x_2)$ and here $\nabla\times=(\partial_1,\partial_2,0)^\trans\times$ and $\nabla\cdot(\partial_1,\partial_2,0)^\trans\cdot\,$.}
$$n=2\quad\mbox{or}\quad n=3,$$
where $\upepsilon_0$ and $\upmu_0$ are the permittivity and the permeability of the free space respectively, $\cE=(\cE_1,\cE_2,\cE_3)$ and $\cH=(\cH_1,\cH_2,\cH_3)$ are the electric and the magnetic field respectively, and the electric displacement field $\cD=(\cD_1(\cE),\cD_2(\cE),\cD_3(\cE))$ depends nonlinearly on the electric field. In the usual case of Kerr isotropic media, one has
\begin{equation}\label{Kerr}
\cD(\cdot,\cE)(x,t)=\epsilon_1(x)\cE(x,t)+\epsilon_3(x)(\cE(x,t)\cdot\cE(x,t))\cE(x,t).
\end{equation}

In the following we set $\upepsilon_0=\upmu_0=1$ without loss of generality. Here we consider two different materials in $\R^n_-$ and $\R^n_+$ respectively, divided by the interface $\Gamma$. We model this configuration by allowing a discontinuity of the linear and cubic susceptibility coefficients $\epsilon_1$ and $\epsilon_3$ across the interface, so that we have
$$\epsilon_j:=\epsilon_j^-\chi_{\R^n_-}+\epsilon_j^+\chi_{\R^n_+}\qquad\mbox{for}\;\;j\in\{1,3\},$$
where $\chi_\Omega$ denotes the characteristic function of the set $\Omega$. We make the usual assumption that $\epsilon_1,\epsilon_3$ are positive and bounded: $0<d\leq\epsilon_1^\pm,\epsilon_3^\pm\in L^\infty(\R^n_\pm)$ for some constant $d$. Because of the inhomogeneity of the material, Maxwell equations \eqref{Max} have to be coupled with suitable interface conditions. In particular, since $\Gamma=\{x_1=0\}$, one has to prescribe
\begin{equation}\label{Max_IFCs}
\llbracket\cE_2\rrbracket=\llbracket\cE_3\rrbracket=0,\quad\;\,\llbracket\cH_1\rrbracket=\llbracket\cH_2\rrbracket=\llbracket\cH_3\rrbracket=0,\quad\;\,\llbracket\cD_1\rrbracket=0,
\end{equation} 
see e.g. \cite[Section 33-3]{Feynman}.

If the material is homogeneous or periodic in the $x_2,\,x_3$ variables and one considers localised quasi-monochromatic signals, it is a standard practice in the physics literature to use an approximation by travelling pulses modulated by an envelope, see e.g. \cite{LS,MS}. A formal ansatz for an asymptotic approximation of a solution $(\cE,\cH)$ of the Maxwell system \eqref{Max} is then
\begin{equation}\label{Uans0}
\cU_{ans}(x,t):=\begin{pmatrix}\cE_{ans}(x,t)\\\cH_{ans}(x,t)\end{pmatrix}:=\varepsilon A(\varepsilon(x_{||}-\nu_gt),\varepsilon^2 t)m_\cU(x_1)\e^{\ri(k_0\cdot x_{||}-\omega_0t)}+\text{c.c.}\footnote{Here c.c. denotes the complex conjugate.},
\end{equation}
where $x_{||}:=(x_2,x_3)$ denotes the components tangential to the interface, $\nu_g\in\R^2$ is the group-velocity of the mode, $k_0\in\R^{n-1}$ is the wave vector, $\omega_0\in\R$ represents a frequency, and $0<\varepsilon\ll1$ is a small asymptotic parameter. The vector function $m_\cU:\R\to\C^6$ is an exponentially decaying ``interface eigenvector'' of a linear problem associated with \eqref{Max} (i.e. $\epsilon_3\equiv0$), whereas the profile $A:\R^{n-1}\times\R\to\C$ has to solve a suitable nonlinear amplitude equation of Schr\"{o}dinger kind in the real slow variables $X=\varepsilon(x_{||}-\nu_gt)$ and $T=\varepsilon^2 t$, see e.g. \cite{KSM,DR,MS}. Therefore, $A$ is smooth and exponentially decaying when sufficient conditions are assumed at the initial time $T=0$.

The above ansatz $\cU_{ans}$ has the property that the divergence equation for the $\cH$-field, i.e. $\nabla\cdot\cH=0$, is automatically satisfied for all times thanks to the choice of $m_\cU$, see Section 5. Moreover, the choice of the profile $A$ is such that when $\cU_{ans}$ is inserted in the Maxwell system, it produces a sufficiently small residuum in $\varepsilon$ for large time scales $t\in[0,T_0(\varepsilon)]$ with $T_0(\varepsilon)\in \cO(\varepsilon^{-2})$ as $\varepsilon \to 0$. Finally, one can notice that 
\begin{equation*}
\partial_t(\nabla\cdot\cD)=\nabla\cdot\partial_t\cD=\nabla\cdot(\nabla\times\cH)=0,
\end{equation*}
which means that whenever the divergence of the $\cD$-field vanishes at $t=0$, it vanishes for all times. The point is that $\nabla\cdot\cD(\cdot,\EA)\not=0$ for 
\begin{equation}\label{Eans0}
\EA(x):=\cE_{ans}(x,0)=\varepsilon A(\varepsilon x_{||},0)m_\cE(x_1)\e^{\ri k_0\cdot x_{||}}+\text{c.c.}\ ,
\end{equation}
where $m_{\cE}$ denotes the vector of the first three components of $m_{\cU}$. Nonetheless one expects that this divergence condition should become satisfied by means of a suitable small correction of the ansatz above. This is precisely the motivation behind the present work. In other words, our aim is to find a vector field $V:\R^n\to\R^n$ which is asymptotically smaller in $\varepsilon$ than $\EA$ and such that
\begin{equation}\label{E0}
\cE^{(0)}:=\EA+V
\end{equation}
solves
\begin{equation}\label{initial_div}
\nabla\cdot\cD(\cdot,\cE^{(0)})=0,
\end{equation}
so that it can be used as an initial condition for the electric field. Of course, $\cE^{(0)}$ has also to fulfil conditions \eqref{Max_IFCs}, i.e. 
\begin{equation}\label{Max_IFCs_korr}
\big\llbracket\cE^{(0)}_2\big\rrbracket=\big\llbracket\cE^{(0)}_3\big\rrbracket=\llbracket\cD_1(\cE^{(0)})\rrbracket=0.
\end{equation}
This correction of higher-order in $\varepsilon$ will not derange the behaviour of $\cE_{ans}$ in the curl-equations. Indeed, if this step is achieved, one then prescribes $\cU_{ans}(x,0)+(V,0)^\trans(x)$ as the initial condition for the electromagnetic field, which then satisfies the two divergence equations in \eqref{Max}, and which is asymptotically close to the initial value of the explicit ansatz $\cU_{ans}$. Having that in hand, one can proceed following the strategy of \cite{SS} and get the existence of an exact solution of the Maxwell equations \eqref{Max} which is well approximated by $\cU_{ans}$ for a large time scale. This part of the analysis is presented in detail in the forthcoming \cite{Daniel}.

\paragraph{General problem and aims} We rewrite our problem \eqref{eq_quasi_gen_D}-\eqref{relationD(E)}-\eqref{IFCs} as
\begin{subnumcases}{}
-\nabla\cdot a(x,\nabla\phi)=b(x) & $\mbox{in}\quad\R^n_\pm,$ \label{eq}
\\
\left\llbracket\left(a(x,\nabla\phi)+\epsilon_1U_0\right)\cdot e_1\right\rrbracket=0, \label{IFCs_eq}
\end{subnumcases}
with $n\geq 2$,
\begin{equation*}
a(x,\nabla\phi):=\epsilon_ff(U_0+\nabla\phi)+\epsilon_1\nabla\phi \quad (= \cD(U_0+\nabla \phi) - \epsilon_1 U_0)
\end{equation*}
and
\begin{equation*}
b:=\nabla\cdot(\epsilon_1U_0),
\end{equation*}
and where $U_0:\R^n\to\R^n$ is a given vector field which satisfies the (linear) transmission condition $\left\llbracket\epsilon_1U_0\cdot e_1\right\rrbracket=0$. Instead of \eqref{IFCs_eq} we often use the shorter notation  $\left\llbracket\cD_1(U_0+\nabla\phi)\right\rrbracket=0$ as in \eqref{IFCs}.

The Maxwell problem described above corresponds then to the case $n=2$ or $n=3$, $U_0=\EA$, $f(v)=|v|^2v$, when we look for an irrotational vector field $V$, i.e. 
\begin{equation*}
V=\nabla\phi
\end{equation*}
for some $\phi:\R^3\to\R$. Note that this formulation does not consider the first two transmission conditions $\big\llbracket\cE^{(0)}_2\big\rrbracket=\big\llbracket\cE^{(0)}_3\big\rrbracket=0$, but they will be automatically satisfied in our application, see Section \ref{Section_application}. 

By \eqref{eq}-\eqref{IFCs_eq} we aim to consider a more general setting of a second-order transmission problem in $\R^n$ which involves a quasilinear operator, namely the divergence of both a linear term (which therefore produces an anisotropic Laplacian) and a nonlinear term, the growth of which is superlinear and behaves like a power. Hence, our difficulties in solving such a problem are threefold: we have to deal with a sort of anisotropic $(p,2)$-Laplacian operator, with the unboundedness of the domain $\R^n_\pm$ and furthermore with transmission conditions on the interface $\Gamma$. As a consequence, we first need to find a suitable functional setting in which to prove the existence of a solution $\phi$. Then, under additional regularity conditions on the coefficients on $\R^n_\pm$, we prove regularity estimates for $\nabla\phi$ which \textit{do not involve} $\|U_0\|_2$. This is of great importance for the above application: $\|U_0\|_2=\|\EA\|_2$ is asymptotically of order $\cO(\varepsilon^{1/2})$ and we are looking for $V=\nabla\phi$ so that it is a correction of $U_0$, i.e. $\|V\|_2=o(\varepsilon^{1/2})$ and for which \eqref{E0}-\eqref{initial_div} hold. This means that the estimates for $\|\nabla\phi\|_2$ (and for its derivatives) should not depend on $\|U_0\|_2$. For details, see Section \ref{Section_application}.

\paragraph{Previous results}
As briefly mentioned, problem \eqref{eq}-\eqref{IFCs_eq} may be classified as a quasilinear transmission problem driven by an operator with a $(p,2)$-growth.

A problem is referred to as \textit{transmission} (or \textit{diffraction}) when the domain in which the equation is defined is split in two or more subdomains in which the coefficients are regular, while at the interfaces they present jump discontinuities, and here the behaviour of the solutions is driven by some compatibility conditions. One usually imposes a condition on the jump of the solutions as well as on the normal derivative with respect to the interface. This class of problems is of great importance for physics and other applied sciences, since they can be derived from many models in which composite materials are treated, not only in the context of propagation of electromagnetic signals, but also e.g. in crack problems \cite{crack}, thermodynamics \cite{MR_thermo} and mathematical biology \cite{BHR,bio}. Moreover, a large amount of papers exist on the numerics of such problems, see e.g. \cite{CZ,NS,CDN,Wang} and references therein.

The theoretical treatment began in the '60s with the classical works of the Russian school by Ole\u{\i}nik, Ladyzhenskaya, Rivkind, Ural'tseva, and coauthors, see e.g. \cite{Oleinik,LRU,LU,LSU}, where linear and quasilinear elliptic and parabolic transmission problems with smooth interfaces were first considered. In these works the main aim was to show that the solutions of the weak formulation are indeed classical, i.e. of class $C^{2,\alpha}$ in the subdomains where the coefficients are regular, and locally of class $C^{1,\alpha}$ near the interface, so that the transmission conditions are satisfied pointwise. Such results were later on improved and extended in \cite{Lorenzi,KL,EFK,TL2}; here, however, the growth of the function $a(x,z,\xi)$ with respect to the gradient-variable $\xi$ is always supposed linear. In \cite{Liu} the authors analyse the case of two adjacent materials which behave according to a model with nonlinear growth of power-type in the gradient-variable $\xi$ (which may also be different from material to material); such results have been later improved in \cite{Knees,KS}. More recently, quasilinear transmission problems even with a wild growth of the function $a$ in the gradient-variable have been investigated in \cite{BBS} by means of Orlicz spaces techniques. We also mention an other direction of research for transmission problems focused on the analysis of non-smooth interfaces \cite{Borsuk}.

All the mentioned results (with the exception of \cite{Lorenzi}) deal with the case of a bounded domain. Equation \eqref{eq}, instead, is set over the whole $\R^n$: this means, roughly speaking, that the growth in the function $a$ with respect to $\xi$ influences the choice of the function space. Having in mind our Maxwell setting, we assume that $a(x,z,\xi)=a(x,\xi)$ and $a(x,\cdot)$ consists of a linear term and power-like nonlinear terms. Therefore the operator we are dealing with can be thought of as a sum of a Laplacian and a $p$-Laplacian. Quasilinear problems of this kind, in the literature called \textit{$(p,q)$-Laplacian} or \textit{double phase} problems, have recently gained a lot of attention in the mathematical community regarding existence and regularity issues. Focusing here on those contributions in which the equation is posed on the whole $\R^n$, we mention \cite{LL,BB,BCS}, in which also a nonlinear right-hand side is considered, and \cite{BC,AdAP,BF}, where in addition a more general double phase operator is considered. In these works the techniques are mainly variational, taking advantage of the mountain pass- or linking-type geometry which the functional associated with the equation enjoys. Orlicz spaces come into play in \cite{BC,AdAP} due to the generality of the operator considered. In \cite{Benci} double phase problems are shown to emerge in the context of static solutions for Lorentz invariant hyperbolic equations on $\R^{3+1}$. On the other hand, the study of \textit{local regularity} for $(p,q)$-Laplacian equations initiated by Marcellini in \cite{Mar} significantly grew in the last two decades, and we refer to \cite{M,MR,Mar2} for an overview on the subject. The focus of these works is on finding sharp conditions on the growth of $a(x,\cdot)$ with respect to the dimension $n$, such that the solutions of \eqref{eq} are locally of class $C^{1,\alpha}$ in the case the coefficients are at least continuous.

We aim to prove standard Sobolev regularity for our solution (which is enough for our application) by means of the well-established method of the difference quotient. Nevertheless we also provide global estimates which are independent on the $L^2$-norm of the given vector field $U_0$ as described above. This complicates the analysis.
\vskip0.2truecm
The rest of the paper is organised as follows. In Section \ref{Section_Var_Ass_Funct} we motivate and describe the suitable functional setting for problem \eqref{eq}-\eqref{IFCs_eq} and we give the precise statements of the assumptions and of our main results. In Section \ref{Section_existence} we prove existence and an estimate for $\nabla\phi$ in terms of the right-hand side $b$ of \eqref{eq} and of $U_0$ but independent of $\|U_0\|_2$. Section \ref{Section_derivatives} contains the proof of the estimates for the derivatives of $\nabla\phi$ in which the discontinuity of the coefficients at the interface comes into play. We shortly illustrate then how analogous estimates are also achievable for higher derivatives of $\nabla\phi$. Finally in Section \ref{Section_application} we apply such results to the original divergence problem \eqref{E0}-\eqref{initial_div}-\eqref{Max_IFCs_korr} arisen in the context of Maxwell equations described in the Introduction and we provide a numerical verification of our results. The short Section \ref{Section_openproblems} contains a discussion about questions which are left open by our analysis, and concludes the paper.

\paragraph{Notation} As usual, for $\Omega\subset\R^n$, $C^\infty_0(\Omega)$ denotes the space of $C^\infty$ functions with compact support in $\Omega$ and $\mathfrak{D}'(\Omega)$ its dual space. $L^p(\R^n)^m$ is the Lebesgue space of $p$-integrable functions over $\R^n$ which take values in $\R^m$. Of course, when $m=1$ we omit the second exponent. The norm in $L^p$-spaces for $1\leq p\leq\infty$ is always denoted by $\|\cdot\|_p$ (i.e. for both scalar and vector case) and $p':=\frac p{p-1}$ and $p^*:=\frac{np}{n-p}$ denote the conjugate exponent and the critical exponent for the Sobolev embedding $W^{1,p}(\R^n)\hookrightarrow L^{p^*}(\R^n)$, respectively. We write $f\in L^p(\R^n_\pm)$ if $f=f_-\chi_{\R^n_-}+f_+\chi_{\R^n_+}$ with $f_\pm\in L^p(\R^n_\pm)$. The spaces $W^{1,p}(\R^n_\pm)$ are analogously defined. The ball with centre $x_0\in\R^n$ and radius $r>0$ is denoted by $B_r(x_0)$. Moreover, we refrain from writing the domain of integration when integration is meant on the whole $\R^n$. The differential $\, dx$ will be often omitted.

\section{Variational structure, assumptions and functional setting}\label{Section_Var_Ass_Funct}

\paragraph{Variational formulation of the problem}
Proceeding formally and postponing to the next paragraph the precise assumptions on $U_0$ and $f$, we start by showing that if we suppose that $f$ is irrotational, i.e. $f=\nabla F$ for some $F\in C^1(\R^n)$, problem \eqref{eq} has a variational structure with the energy functional
\begin{equation*}
\cJ(\phi)=\int\epsilon_fF(U_0+\nabla\phi)+\int\epsilon_1\frac{|\nabla\phi|^2}2+\int\epsilon_1U_0\!\cdot\!\nabla\phi.
\end{equation*}
Notice that, due to the particular form of the right-hand side $b$, $\cJ$ involves only the gradient of $\phi$ and not the function $\phi$ itself. Indeed, we claim that the Euler-Lagrange equation associated with $\cJ$ is \eqref{eq} and that critical points of $\cJ$ satisfy the interface condition \eqref{IFCs_eq}. For any $\eta\in C^\infty_0(\R^n)$ we get
\begin{equation*}
\cJ'[\phi](\eta)=\int\epsilon_ff(U_0+\nabla\phi)\ccdot\nabla\eta+\int\epsilon_1(U_0+\nabla\phi)\ccdot\nabla\eta.
\end{equation*}
Splitting each integral over the two half-spaces $\R^n_\pm$ and applying formally Gau\ss{}'s theorem to each of them, we get
\begin{equation*}
\begin{split}
\int\epsilon_1(U_0+\nabla\phi)\ccdot\nabla\eta&=\int_{\R^n_-}\epsilon_1^-(U_0+\nabla\phi)\ccdot\nabla\eta+\int_{\R^n_+}\epsilon_1^+(U_0+\nabla\phi)\ccdot\nabla\eta\\
&=-\int\nabla\cdot(\epsilon_1(U_0+\nabla\phi))\eta-\int_\Gamma\big((\epsilon_1(U_0+\nabla\phi))_+-(\epsilon_1(U_0+\nabla\phi))_-\big)\cdot e_1\eta,
\end{split}
\end{equation*}
where $u_\pm(x):=\lim_{t\to0^\pm}u(x + t e_1)$, $x \in \Gamma$. Analogous computations hold for the first term
\begin{equation*}
\begin{split}
\int\epsilon_ff(U_0+\nabla\phi)\ccdot\nabla\eta&=-\int\nabla\cdot\left(\epsilon_ff(U_0+\nabla\phi)\right)\eta\\
&\quad-\int_\Gamma\left(\left(\epsilon_ff(U_0+\nabla\phi)\right)_+-\left(\epsilon_ff(U_0+\nabla\phi)\right)_-\right)\cdot e_1\eta.
\end{split}
\end{equation*}
All in all smooth critical points of $\cJ$ satisfy
\begin{equation*}
\int\nabla\cdot\left(\epsilon_ff(U_0+\nabla\phi)+\epsilon_1\nabla\phi\right)\eta\,+\int\nabla\cdot(\epsilon_1U_0)\eta\,+\int_\Gamma\!\!\Big(\cD(\cdot,U_0+\nabla\phi)_+-\cD(\cdot,U_0+\nabla\phi)_-\Big)\,\cdot\,e_1\eta=0.
\end{equation*}
Taking now $\eta\in C^\infty_0(\R^n)$ with $\supp\,\eta\subset\R^n_\pm$, the boundary integral vanishes and we obtain that the critical points of $\cJ$ satisfy \eqref{eq} in $\R^n_\pm$. If, on the other hand, $\supp\,\eta\cap\Gamma\neq\emptyset$, then the boundary integral tells us that the quantity $\big(\cD(\cdot,U_0+\nabla\phi)_+-\cD(\cdot,U_0+\nabla\phi)_-\big)\cdot e_1=\llbracket \cD_1(\cdot,U_0+\nabla\phi)\rrbracket$ vanishes a.e., which is condition \eqref{IFCs_eq}.

As we work in a Sobolev space, we need to interpret the interface condition via the \textit{normal trace} of $\cD(\cdot,U_0+\nabla\phi)$ on $\R^n_+$ and $\R^n_-$. The standard definition of the normal trace for $\Hdiv$-functions is not suitable in our setting because our solutions do not satisfy $\cD(\cdot,U_0+\nabla \phi)\in L^2(\R^n)$. Instead we define the trace in the sense of distributions.
\begin{defn}\label{def_trace}
	For
	$$S_\pm:=\left\{\cD\in L^1_{loc}(\R^n_\pm)^n\,\bigg|\, \nabla\cdot \cD\in  L^1_{loc}(\R^n_\pm)\right\}$$
	we define the traces
	\begin{equation*}
		\cT_\pm: S_\pm\to\mathfrak{D}'(\Gamma),\quad \big(\cT_\pm\cD\big)[\psi]:=\mp\int_{\R^n_\pm}\nabla\cdot\cD\,\hat{\psi} + \cD\cdot \nabla \hat{\psi}\quad\mbox{for all}\;\,\psi\in C^\infty_0(\Gamma),
	\end{equation*}
where $\hat{\psi}\in C^\infty_0(\overline{\R^n_\pm})$ such that $\hat{\psi}|_\Gamma = \psi$.
\end{defn}
The existence of $\hat{\psi}$ was proved, e.g. in \cite{Whitney}. The fact that the trace is independent of the choice of $\hat{\psi}$ follows directly from the definition of the weak divergence. With this definition the interface condition \eqref{IFCs_eq} in trace sense is
\begin{equation}\label{T+=T-}
\cT_+\cD(U_0+\nabla\phi)=\cT_-\cD(U_0+\nabla\phi) \quad \text{in} \ \mathfrak{D}'(\Gamma).
\end{equation}
As the above calculation shows, if $\nabla\cdot\cD(\cdot,U_0+\nabla\phi)$ exists in the classical sense, then \eqref{T+=T-} is equivalent to \eqref{IFCs_eq} pointwise on $\Gamma$.

\paragraph{Assumptions}
Trying to retain the essential features of the Maxwell context described in the introduction, we assume the following conditions for $f$, $\epsilon_1$, $\epsilon_f$, $U_0$ throughout the paper.
\begin{enumerate}
	\item[(H0)] $\epsilon_1,\epsilon_f\in L^\infty(\R^n)$ and there exists a constant $d$ such that $\epsilon_1(x),\epsilon_f(x)\geq d>0$ a.e. in $\R^n$;
	\item[(F0)] there exists $F\in C^1(\R^n)$ convex so that $F(0)=0$, $f=\nabla F$, and $F(v)\geq\mu_p|v|^{p+1}$ for some $\mu_p>0$;
	\item[(F1)] there exist exponents $1<\alpha\leq p$ and constants $0<\lambda_p\leq\Lambda_p$ and $\Lambda_\alpha\geq 0$ such that for all $v\in \R^n$
	\begin{enumerate}[i)]
		\item $|f(v)|\leq\Lp|v|^p+\La|v|^\alpha\,$;
		\item $f(v)\ccdot v>\lp|v|^{p+1}$;
	\end{enumerate}
	\item[(A0)] $U_0\in L^2(\R^n)^n\cap L^{p+1}(\R^n)^n$;
	\item[(A1)] $b:=\nabla\cdot(\epsilon_1U_0)\in L^2(\R^2)\cap L^1(\log,\R^2)$ if $n=2$ and $b\in L^2(\R^n)\cap L^{\frac{2n}{n+2}}(\R^n)$ if $n\geq 3$, where $L^1(\log,\R^2):=\left\{\varphi\in L^1_{loc}(\R^2)\,\big|\,\|\varphi\|_{L^1(\log,\R^2)}:=\int\log(2+|x|)|\varphi(x)|\,dx<\infty\right\}$. To unify the two cases, we write 
\begin{equation}\label{OS_norm}
\iii b:=\begin{cases}
\|b\|_{L^1(\log,\R^2)}&\mbox{if}\;\,n=2,\\
\|b\|_{\frac{2n}{n+2}}&\mbox{if}\;\,n\geq3;
\end{cases}
\end{equation}
	\item[(A2)] $U_0$ does not weakly solve the nonlinear equation $\nabla\cdot\cD(\cdot,\cE)=0$, where $\cD(\cdot,\cE)$ is defined in \eqref{relationD(E)} (i.e. $\phi=0$ is not a solution in the sense of Definition \ref{weak_sol}), but satisfies the transmission condition $\left\llbracket\epsilon_1U_0\cdot e_1\right\rrbracket=0$ (in the sense of Definition \ref{def_trace}).
\end{enumerate}

To obtain the existence and $L^2$-estimates for the derivatives of $\nabla\phi$, we need moreover
\begin{enumerate}
	\item[(H1)] $\epsilon_1,\epsilon_f\in W^{1,\infty}(\R^n_\pm)$;
\end{enumerate}
\begin{enumerate}
	\item[(F2)] $f\in C^1(\R^n)^n$ and there exist constants $0<\tlp\leq\tLp$ and $\tLa\geq 0$ such that the Jacobi matrix $J_f$ of $f$ satisfies for all $v,\xi\in\R^n$
	\begin{enumerate}[i)]
		\item $|J_f(v)|\leq\tLp|v|^{p-1}+\tLa|v|^{\alpha-1}$;
		\item $\left(J_f(\xi)v\right)\ccdot v\geq\tlp|\xi|^{p-1}|v|^2\,$;
	\end{enumerate}
	\item[(A3)] $\partial_kU_0\in L^2(\R^n)^{n\times n}\cap L^{p+1}(\R^n)^{n\times n}$ for $k\in\{2,\dots,n\}$, and $\partial_1U_0\in L^2(\R^n_\pm)^{n\times n}\cap L^{p+1}(\R^n_\pm)^{n\times n}$;
	\item[(A4)] $\partial_k b\in L^2(\R^n)$, $k\in\{2,\dots,n\}$, and $\partial_1b\in L^2(\R^n_\pm)$.
\end{enumerate}

\begin{remark}
	Assumptions (H0)-(H1) on the coefficients $\epsilon_1,\epsilon_f$ give us the structure of a transmission problem.
\end{remark}

\begin{remark}
	The simplest (and physically relevant) model nonlinearity is $f(v)=|v|^{p-1}v$ with $p>1$, but also nonlinearities of the kind $f(v)=\sum_{i=1}^N|v|^{q_i}v$ with $q_i>0$ satisfy (F0)-(F2) with $\alpha=1+\min_i q_i$ and $p=1+\max_i q_i$ and are allowed.
\end{remark}
\begin{remark}
	Assumption (A2) reflects the situation that we encounter in our application to the Maxwell setting, where the vector field $U_0$, being just a solution of the linearised equation and not of the nonlinear equation $\nabla\cdot\cD(\cdot,\cE)=0$, still satisfies the transmission condition $\llbracket\cD_1(\cdot, U_0)\rrbracket=0$. Notice that such a transmission condition allows to apply Gau\ss{}'s theorem and transform the third term in the functional $\cJ$ into $\int b\phi$, which exists thanks to assumption (A1), see Section 3.
\end{remark}

\paragraph{The main results} From assumptions (H0), (F1), and (A0), it is easy to see that the suitable functional space in which to look for critical points of the functional $\cJ$ is
\begin{equation}\label{D24}
\mathfrak D_{2,p+1}:=D^{1,2}_0(\R^n)\cap D^{1,p+1}_0(\R^n),
\end{equation}
where for $q\geq 1$
\begin{equation*}
D^{1,q}_0(\R^n):=\overline{C^\infty_0(\R^n)}^{\,|\cdot|_{1,q}}\qquad\mbox{with the norm}\qquad|u|_{1,q}:=\|\nabla u\|_q
\end{equation*}
is the \textit{homogeneous Sobolev space}, sometimes in the literature referred to as the \textit{Beppo Levi space}. The norm on $\mathfrak D_{2,p+1}$ is defined as $\|\cdot\|_\fD:=|\cdot|_{1,2}+|\cdot|_{1,p+1}$. 

\begin{defn}\label{weak_sol}
	We say that $\phi\in\fD_{2,p+1}$ is a \textit{weak solution} of problem \eqref{eq}-\eqref{IFCs_eq} if
	\begin{equation}\label{wf}
	\int\epsilon_ff(U_0+\nabla\phi)\ccdot\nabla\eta+\int\epsilon_1(U_0+\nabla\phi)\ccdot\nabla\eta=0\qquad\mbox{for all}\,\;\eta\in\fD_{2,p+1}.
	\end{equation}
\end{defn}

Weak solutions indeed satisfy the interface condition \eqref{IFCs_eq} in the trace sense, which follows from the definition and from the fact that $\nabla \cdot \cD(U_0+\nabla\phi)=0$ pointwise almost everywhere, see Theorem \ref{Thm_der_phi} and its proof.

Notice that \eqref{eq} depends just on the gradient of the unknown function, and not on the function itself. This is the first difference between our analysis and the quasilinear problems on the whole $\R^n$ of the same kind which appear in the literature: indeed the $p$-Laplace-kind operator is often coupled with the term $|u|^{p-2}u$, so that one can work within the functional framework of classical Sobolev spaces, see e.g. \cite{LL,BF,BCS}.
\vskip0.2truecm
Our main results are the following:
\begin{thm}\label{Thm_phi}
	Let $\epsilon_1,\,\epsilon_f,\,f,\,U_0$ satisfy assumptions (H0), (F0)-(F1), (A0)-(A2). Then there exists a nontrivial minimum $\phi$ of the functional $\cJ$ in $\fD_{2,p+1}$ and there holds
	\begin{equation}\label{estimate_phi}
	\int|\nabla\phi|^2+\int|U_0+\nabla\phi|^{p-1}|\nabla\phi|^2\leq C\left(\|U_0\|_{p+1}^{p+1}+\|U_0\|_{\alpha+1}^{\alpha+1}+\|b\|_2^2+\iii b^2\right),
	\end{equation}
	where the constant $C$ depends only on $\Lp$, $\lp$, $\La$, $d$, $\|\epsilon_1\|_\infty$, and $\|\epsilon_f\|_\infty$. Moreover, $\phi$ is a weak solution of \eqref{eq}-\eqref{IFCs_eq}.
\end{thm}

\begin{thm}\label{Thm_der_phi}
	Let $\epsilon_1,\,\epsilon_f,\,f,\,U_0$ satisfy assumptions (H0)-(H1), (F0)-(F2), (A0)-(A4) and let $\phi$ be a minimiser of Theorem \ref{Thm_phi}. Then for $\nabla\phi$ the tangential (with respect to $\Gamma$) derivatives $\partial_k\nabla\phi$, $k\in\{2,\dots,n\}$, satisfy $\partial_k\nabla\phi\in L^2(\R^n)^n$ and
	\begin{equation}\label{estimate_der_phi_tg}
	\begin{split}
	\int|\partial_k\nabla\phi|^2+\int|U_0+\nabla\phi|^{p-1}|\partial_k\nabla\phi|^2&\leq C_k\Big(\|U_0\|_{p+1}^{p+1}+\|\partial_kU_0\|_{p+1}^{p+1}+\|U_0\|_{\alpha+1}^{\alpha+1}\\
	&\quad+\|\partial_kU_0\|_{\alpha+1}^{\alpha+1}+\|b\|_2^2+\iii b^2+\|\partial_kb\|_2^2\Big).
	\end{split}
	\end{equation}
	The normal derivative $\partial_{11}\phi$ exists in $L^2(\R^n_\pm)$ and there holds
	\begin{equation}\label{estimate_der_phi_n}
	\begin{split}
	\int_{\R^n_\pm}|\partial_{11}\phi|^2+\int_{\R^n_\pm}|U_0+\nabla\phi|&^{p-1}|\partial_{11}\phi|^2\leq C_1\Big(\|U_0\|_{p+1}^{p+1}+\sum_{k=1}^n\|\partial_kU_0\|_{L^{p+1}(\R^n_\pm)}^{p+1}+\|U_0\|_{\alpha+1}^{\alpha+1}\\
	&\quad+\sum_{k=1}^n\|\partial_kU_0\|_{L^{\alpha+1}(\R^n_\pm)}^{\alpha+1}+\|b\|_2^2+\iii b^2+\sum_{k=2}^n\|\partial_kb\|_{L^2(\R^n_\pm)}^2\Big).
	\end{split}
	\end{equation}
	The constants $C_k,C_1$ depend only on $\Lp$, $\lp$, $\La$, $\tLp$, $\tlp$, $\tLa$, $d$, $\|\epsilon_1\|_{W^{1,\infty}(\R^n\setminus \Gamma)}$, and $\|\epsilon_f\|_{W^{1,\infty}(\R^n\setminus \Gamma)}$.
	
	Moreover, equation \eqref{eq} is satisfied pointwise almost everywhere and the interface condition \eqref{IFCs_eq} holds in the trace sense of Definition \ref{def_trace}.
\end{thm}
\begin{remark}\label{Rmk_norm_der}
In estimate \eqref{estimate_der_phi_n} only the first component $\partial_{11}\phi$ of the normal derivative $\partial_1\nabla\phi$ appears. Indeed, for the remaining components there holds $\partial_{1k}\phi=\partial_{k1}\phi$ a.e. in $\R^n_\pm$, $k\in\{2,\dots,N\}$, and thus estimate \eqref{estimate_der_phi_tg} can be directly applied. Hence, only $\partial_{11}\phi$ has to be studied separately. For details, see Section \ref{Sec_norm_der}.
\end{remark}
\begin{remark}
	An upper bound for $\|\nabla\phi\|_2^2+\|\nabla\phi\|_{p+1}^{p+1}$ analogous to \eqref{estimate_phi} with the same terms on the right-hand side can be obtained by \eqref{nablaphi_p+1_pt2}. However, we preferred to state Theorem \ref{Thm_phi} this way to have a better comparison with the left-hand side of the inequalities \eqref{estimate_der_phi_tg}-\eqref{estimate_der_phi_n}.
\end{remark}

\begin{remark}\label{Remark_smallness}
	Notice that the quantity on the right-hand side of \eqref{estimate_phi} does \textit{not} involve $\|U_0\|_2$. This will be crucial in our application to the asymptotic Maxwell problem in Section \ref{Section_application} in order to show that $\nabla\phi$ is actually an $L^2$-correction of $U_0$, as desired.
\end{remark}

\begin{remark}\label{Remark_higher_regularity}
	If one prescribes higher regularity on $\epsilon_1,\,\epsilon_f,\,U_0,\,b$, as well as control on the higher derivatives of $f$ analogously to assumptions (F1)-(F2), one may further infer similar estimates for the higher derivatives of $\nabla\phi$ to the ones in \eqref{estimate_der_phi_tg}-\eqref{estimate_der_phi_n}. The proof is then analogous to that of Theorem \ref{Thm_der_phi}.
\end{remark}

\paragraph{The functional setting: homogeneous Sobolev spaces}
Before giving the proofs of Theorems \ref{Thm_phi}-\ref{Thm_der_phi} in the next sections, for the sake of a self-contained exposition and in order to state some properties of such spaces that are needed in the sequel, we take now a small detour and define homogeneous Sobolev spaces following \cite[Chapter II]{Galdi}, see also \cite[Sec 6]{Kuro}. For $\Omega\subset\R^n$ arbitrary measurable domain, $m\in\N$, and $q\geq1$, we start by defining
$$D^{m,q}(\Omega):=\{u\in L^1_{loc}(\Omega)\,|\,D^\ell u\in L^p(\Omega),\, |\ell|=m\}$$
with the norm $|u|_{m,q}:=\|\nabla^mu\|_q$. If $\Omega$ is unbounded, a control on the norm of the highest derivative does \textit{not} imply a control on the norms of the function itself and of all other derivatives. In other words $D^{m,q}(\Omega)\setminus W^{m,q}(\Omega)\neq\emptyset$. Nevertheless, if $q>1$ any function $u\in D^{m,q}(\Omega)$ belongs to $W^{m,q}_{loc}(\Omega)$ and
\begin{equation}\label{estimateLpD1p}
\|u\|_{ W^{m,q}(\Omega')}\leq C\left(\sum_{|\ell|=m}\|D^\ell u\|_{L^q(\Omega')}+\|u\|_{L^1(\Omega')}\right)
\end{equation}
for any $\Omega'\subset\subset\Omega$, see \cite[Lemma II.6.1]{Galdi}.

An immediate disadvantage of these spaces is that, if $P\in\cP^m:=\{\mbox{polynomials of degree}\,\leq m-1\}$, then it is clear that $|u|_{m,q}=|u+P|_{m,q}$. One can avoid this by defining
$$\Do^{m,q}(\Omega):=\{[u]_m\,|\,u\in D^{m,q}(\Omega)\},$$
where $[u]_m:=\{w\in D^{m,q}(\Omega)\,|\,w=u+P,\,\mbox{for some}\,P\in\cP^m\}$ and one can show that $\Do^{m,q}(\Omega)$ is a Banach space. Since $C^\infty_0(\Omega)\subset\Do^{m,q}(\Omega)$ via the natural inclusion $i:u\mapsto i(u):=[u]_m$, the space $D^{m,q}_0(\Omega):=\overline{C^\infty_0(\Omega)}^{|\cdot|_{m,q}}$ is isomorphic to a closed subset of $\Do^{m,q}(\Omega)$, so it is a Banach space too.
\begin{remark} 
	Notice that this means that if we consider $u\in D^{1,q}_0(\R^n)$, then it coincides with an element of $D^{1,q}(\R^n)$, so $u\in L^1_{loc}(\R^n)$. This in turns implies that $u\in W^{1,q}_{loc}(\R^n)$ by \eqref{estimateLpD1p}.
\end{remark}

The spaces $\Do^{m,q}(\Omega)$ and $D^{m,q}_0(\Omega)$ turn out to be separable for $1\leq q<+\infty$ and reflexive for $1<q<\infty$, see \cite[Exercise II.6.2]{Galdi} or \cite[Theorem 2.2]{S}. This implies that our space $\fD_{2,p+1}$ defined in \eqref{D24} is a reflexive Banach space. Indeed, the spaces $D^{1,2}_0(\R^n)$ and $D^{1,p+1}_0(\R^n)$ form a conjugate couple of Banach spaces (i.e. their intersection is dense in both spaces) and therefore, see \cite[Theorem 8.III]{AG65},
\begin{equation*}
\begin{split}
\fD_{2,p+1}''&=\left(D^{1,2}_0(\R^n)\cap D^{1,p+1}_0(\R^n)\right)''=\left(D^{1,2}_0(\R^n)'+ D^{1,p+1}_0(\R^n)'\right)'\\
&=D^{1,2}_0(\R^n)''\cap D^{1,p+1}_0(\R^n)''=D^{1,2}_0(\R^n)\cap D^{1,p+1}_0(\R^n)=\fD_{2,p+1}.
\end{split}
\end{equation*}

Finally, we recall a result which allows us to split a function which belongs to $\Do^{1,q}(\Omega)$ into a ``Sobolev'' term in $W^{1,q}$ and a ``regular'' term which might not be in $W^{1,q}$ but over which we may have a control in some Lebesgue space.
\begin{prop}[\cite{OS}, Theorem 2.2 (iii)]\label{OS_thm}
	For $q\in[1,\infty)$ there exist linear maps $T_0:\Do^{1,q}(\R^n)\to W^{1,q}(\R^n)$ and $T_\infty:\Do^{1,q}(\R^n)\to C^\infty(\R^n)$ such that $[u]=[T_0(u)]+[T_\infty(u)]$ for all $[u]\in\Do^{1,q}$ and
	\begin{equation*}
	\|T_0(u)\|_{W^{1,q}}\leq C|u|_{1,q},\quad|T_\infty(u)|_{1,q}\leq|u|_{1,q},\quad\|\nabla T_\infty(u)\|_\infty\leq|u|_{1,q},
	\end{equation*}
	where $C=C(n)>0$. Moreover $T_\infty$ may be chosen to satisfy
	\begin{equation*}
	\begin{aligned}
	&|T_\infty(u)|(x)\leq C|u|_{1,q}\log(2+|x|),\quad x\in\R^n\quad&\text{if } \; q=n,\\
	&\|T_\infty(u)\|_{q^*}\leq C|u|_{1,q} &\text{if } \; q<n.
	\end{aligned}
	\end{equation*}
\end{prop}

\begin{remark}
	This result motivates our definition of the norm $\iii\cdot$.
\end{remark}

\section{Existence and estimates for \texorpdfstring{$\boldsymbol{\nabla\phi}$}{nabla phi} (proof of Theorem \ref{Thm_phi})}\label{Section_existence}

\paragraph{Notation} In the sequel we often use the symbol $\lesssim$ to indicate that an inequality holds up to a multiplicative constant depending only on the structural constants $n$, $\Lp$, $\lp$, $\La$, $\tLp$, $\tlp$, $\tLa$, $d$, $\|\epsilon_1\|_{W^{1,\infty}(\R^n\setminus \Gamma)}$, and $\|\epsilon_f\|_{W^{1,\infty}(\R^n\setminus \Gamma)}$ but not on $U_0$, $b$, and $\phi$. Moreover $\delta$ is a small positive parameter and $C_\delta$ denotes a positive constant which smoothly depends on $\delta$ and whose value may vary from line to line.

\paragraph{Minimisation}
Under assumptions (H0), (F0)-(F1), and (A0), the existence part of Theorem \ref{Thm_phi} essentially follows by the standard minimisation method.
First, note that the functional $\cJ$ is well-defined on $\mathfrak D_{2,p+1}$, since
\begin{equation*}
\begin{split}
	\left|\int\epsilon_fF(U_0+\nabla\phi)\right|&\leq\|\epsilon_f\|_\infty\int|F(U_0+\nabla\phi)|\\
	&\leq\|\epsilon_f\|_\infty\left(\Lp\int|U_0+\nabla\phi|^{p+1}+\La\int|U_0+\nabla\phi|^{\alpha+1}\right)\\
	&\les \|U_0\|_2^2+\|U_0\|_{p+1}^{p+1}+\|\nabla\phi\|_2^2+\|\nabla\phi\|_{p+1}^{p+1}
\end{split}
\end{equation*}
since $F(0)=0$ and $\alpha\in(1,p]$. Moreover, $\cJ$ is coercive on $\fD_{2,p+1}$. Indeed, by (F0)
\begin{equation*}
\begin{split}
\cJ(\phi)&\geq d\mu_p\|U_0+\nabla\phi\|_{p+1}^{p+1}+\frac d2\|\nabla\phi\|_2^2-\|\epsilon_1\|_\infty\|U_0\|_2\|\nabla\phi\|_2\\
&\geq C\|\nabla\phi\|_{p+1}^{p+1}+\left(\tfrac d2-\delta\right)\|\nabla\phi\|_2^2-C\|U_0\|_{p+1}^{p+1}-C_\delta\|U_0\|_2^2,
\end{split}
\end{equation*}
and the coercivity of $\cJ$ follows by choosing $\delta$ sufficiently small. Finally, the weakly lower semicontinuity is ensured by the convexity of $\cJ$, see e.g. \cite[Theorem 4.5]{Giusti}. Indeed, $\cJ$ is a sum of a squared weighted $L^2$-norm and a nonlinear term driven by $F$, which is convex due to (F0). Since $\fD_{2,p+1}$ is a reflexive Banach space as shown in Section \ref{Section_Var_Ass_Funct}, the direct method of the calculus of variations provides the existence of a global minimum of $\cJ$, which is in particular a weak solution of \eqref{eq}-\eqref{IFCs_eq}. Note that such a minimum is not trivial because otherwise $\nabla\cdot\cD(\cdot,U_0)=0$ would weakly hold and thus violate assumption (A2).

\paragraph{Estimate \eqref{estimate_phi}}
Testing \eqref{wf} with $\eta=\phi$, we find 
\begin{equation}\label{eq_L2f}
\int\epsilon_ff(U_0+\nabla\phi)\ccdot(U_0+\nabla\phi)+\int\epsilon_1|\nabla\phi|^2=\int\epsilon_ff(U_0+\nabla\phi)\ccdot U_0-\int\epsilon_1U_0\ccdot\nabla\phi.
\end{equation}
Since $\epsilon_1$ and $\epsilon_f$ are bounded from below by a positive constant, see (H0), the growth condition (F1.ii) implies that 
\begin{equation}\label{eq_L2f_LHS}
\int\epsilon_ff(U_0+\nabla\phi)\ccdot(U_0+\nabla\phi)+\int\epsilon_1|\nabla\phi|^2\geq d\lp\int|U_0+\nabla\phi|^{p+1}+d\int|\nabla\phi|^2.
\end{equation}
On the other hand, by (F1.i), H\"older inequality with exponents $\frac{p+1}p$ and $p+1$, and a $\delta$-Young inequality we have 
\begin{equation*}
\begin{split}
\left|\int\epsilon_ff(U_0+\nabla\phi)\ccdot U_0\right|&\leq\|\epsilon_f\|_\infty\left(\Lp\int|U_0+\nabla\phi|^p|U_0|+\La\int|U_0+\nabla\phi|^\alpha|U_0|\right)\\
&\leq\delta\int|U_0+\nabla\phi|^{p+1}+C_\delta\int|U_0|^{p+1}+\delta\int|U_0+\nabla\phi|^{\alpha+1}+C_\delta\int|U_0|^{\alpha+1}
\end{split}
\end{equation*}
and, using $x^{\alpha+1}\les x^{p+1}+x^2$ for $x\geq0$, since $\alpha\in(1,p]$, 
\begin{equation*}
\begin{split}
\int|U_0+\nabla\phi|^{\alpha+1}&\les\int|\nabla\phi|^{p+1}+\int|\nabla\phi|^2+\int|U_0|^{\alpha+1}\\
&\les\int|U_0+\nabla\phi|^{p+1}+\int|\nabla\phi|^2+\|U_0\|_{\alpha+1}^{\alpha+1}+\|U_0\|_{p+1}^{p+1}.
\end{split}
\end{equation*}
Therefore
\begin{equation}\label{eq_L2f_RHS}
\begin{split}
\left|\int\epsilon_ff(U_0+\nabla\phi)\ccdot U_0\right|\les\delta\int|U_0+\nabla\phi|^{p+1}+\delta\int|\nabla\phi|^2+C_\delta\|U_0\|_{\alpha+1}^{\alpha+1}+C_\delta\|U_0\|_{p+1}^{p+1},
\end{split}
\end{equation}
where the first two terms on the right-hand side can be absorbed by the quantities on the right of \eqref{eq_L2f_LHS}. It remains to estimate the last term in \eqref{eq_L2f}. Since we aim to obtain an estimate of the $L^2$- and $L^{p+1}$-norms of $\nabla\phi$ independently of the $L^2$-norm of $U_0$, we need to rely on the term $b=\nabla\cdot(\epsilon_1 U_0)$. By (A1)-(A2) we deduce $b\in L^2(\R^n)$ as well as the transmission condition $\left\llbracket\epsilon_1U_0\cdot e_1\right\rrbracket=0$. Hence
\begin{equation}\label{b_Gauss}
-\int\epsilon_1U_0\ccdot\nabla\phi=\int b\phi
\end{equation}
and we aim to show that the right-hand side is well-defined in $L^1$, which is not trivial since we are working in homogeneous Sobolev spaces. To this purpose, we decompose $\phi=u+v$ according to Proposition \ref{OS_thm} with $u:=T_\infty(\phi)$ and $v:=T_0(\phi)$. Hence we have $v\in H^1(\R^n)$ with $\|v\|_{H^1}\leq C|\phi|_{1,2}$ and either $|u(x)|\leq C\log(2+|x|)|\phi|_{1,2}$ if $n=2$ or $\|u\|_{2^*}\leq C|\phi|_{1,2}$ if $n\geq3$. Then
\begin{equation}\label{eq_L2_b}
\begin{split}
\left|\int b\phi\right|&\les\|b\|_2\|v\|_2+\iii b|\phi|_{1,2}\les \left(\|b\|_2+\iii b\right)|\phi|_{1,2}\\
&\les\delta\int|\nabla\phi|^2+ C_\delta\left(\|b\|_2^2+\iii b^2\right).
\end{split}
\end{equation}
Hence, \eqref{eq_L2f}-\eqref{eq_L2_b} yield
\begin{equation*}
\begin{split}
(1-C_1\delta)\int|U_0+\nabla\phi|^{p+1}+(1-C_2\delta)\int|\nabla\phi|^2&\les C_\delta\left(\|U_0\|_{\alpha+1}^{\alpha+1}+\|U_0\|_{p+1}^{p+1}+\|b\|_2^2+\iii b^2\right).
\end{split}
\end{equation*}
for suitable constants $C_1,C_2>0$. To infer the desired estimate \eqref{estimate_phi} it is now sufficient to choose $\delta$ sufficiently small after having observed that
\begin{equation}\label{weighted_term_Lp+1}
\begin{split}
\int|U_0+\nabla\phi|^{p-1}|\nabla\phi|^2&\les\int|U_0+\nabla\phi|^{p+1}+\int|U_0+\nabla\phi|^{p-1}|U_0|^2\\
&\les\int|U_0+\nabla\phi|^{p+1}+\|U_0\|_{p+1}^{p+1},
\end{split}
\end{equation}
where again H\"older inequality with conjugate exponents $\frac{p+1}{p-1}$ and $\frac {p+1}2$ and a $\delta$-Young inequality have been used. This concludes the proof of Theorem \ref{Thm_phi}.
\vskip0.2truecm
For later use, we note here that the ``almost converse'' inequality
\begin{equation}\label{weighted_term_Lp+1_converse}
\int|U_0+\nabla\phi|^{p+1}\les\int|U_0+\nabla\phi|^{p-1}|\nabla\phi|^2+\|U_0\|_{p+1}^{p+1}
\end{equation}
holds. Indeed, by similar arguments as for \eqref{weighted_term_Lp+1} we have
\begin{equation}\label{weighted_term_Lp+1_converse_Schritt1}
\begin{split}
\int|U_0+\nabla\phi|^{p+1}&\leq\int|U_0+\nabla\phi|^{p-1}|U_0|^2+\int|U_0+\nabla\phi|^{p-1}|\nabla\phi|^2\\
&\les\|U_0\|_{p+1}^{p+1}+\int|\nabla\phi|^{p-1}|U_0|^2+\int|U_0+\nabla\phi|^{p-1}|\nabla\phi|^2\\
&\leq2\|U_0\|_{p+1}^{p+1}+\|\nabla\phi\|_{p+1}^{p+1}+\int|U_0+\nabla\phi|^{p-1}|\nabla\phi|^2.
\end{split}
\end{equation}
Now we further estimate
\begin{equation*}
\begin{split}
\|\nabla\phi\|_{p+1}^{p+1}&\les\int|U_0+\nabla\phi|^{p-1}|\nabla\phi|^2+\int|U_0|^{p-1}|\nabla\phi|^2\\
&\les\int|U_0+\nabla\phi|^{p-1}|\nabla\phi|^2+\delta\|\nabla\phi\|_{p+1}^{p+1}+C_\delta\|U_0\|_{p+1}^{p+1}
\end{split}
\end{equation*}
by the $\delta$-Young inequality, so choosing $\delta=\tfrac12$ one gets
\begin{equation}\label{nablaphi_p+1_pt2}
\begin{split}
\|\nabla\phi\|_{p+1}^{p+1}\les\int|U_0+\nabla\phi|^{p-1}|\nabla\phi|^2+\|U_0\|_{p+1}^{p+1}.
\end{split}
\end{equation}
Estimate \eqref{weighted_term_Lp+1_converse} follows by combining \eqref{weighted_term_Lp+1_converse_Schritt1} with \eqref{nablaphi_p+1_pt2}.

\section{Estimates for the derivatives of \texorpdfstring{$\boldsymbol{\nabla\phi}$}{nabla phi} (proof of Theorem \ref{Thm_der_phi})}\label{Section_derivatives}
In the previous section we proved the existence of $\phi$ in $\fD_{2,p+1}$, together with estimate \eqref{estimate_phi}, under the sole assumption that the coefficients $\epsilon_1,\,\epsilon_f$ are bounded from above and below by positive constants. Their regularity was not involved. However, for the existence of the derivatives of $\nabla\phi$ as well as analogous estimates on them, the regularity of $\epsilon_1,\,\epsilon_f$ comes into play. Assumption (H1) guarantees that their weak derivatives exist in the half-spaces $\R^n_\pm$ but they might have a jump at the interface, which is the typical situation in applications. We show that the derivatives of $\nabla\phi$ which are tangential to the interface exist on the whole $\R^n$, while the normal derivative $\partial_{11}\phi$ is well-defined just in the two half-spaces. In both cases we obtain estimates analogous to the one found in the previous section.

We make use of the well-known method of difference quotients. In the context of transmission problems it was employed in \cite{LU} to obtain local regularity for linear equations and then extended in \cite{TL2} to quasilinear problems of the kind \eqref{eq} with a function $a(x,\cdot)$ with at most linear growth. On the other hand, this method was also applied to double phase quasilinear problems in \cite{Mar}. Here we want to merge these two features and obtain analogous estimates for our double phase quasilinear transmission problem \eqref{eq}-\eqref{IFCs_eq}.
\vskip0.2truecm
We start by defining the difference quotient of a function $u:\R^n\to\R$ in the direction $e_k$, $k\in\{1,\dots,n\}$ as
\begin{equation*}
\partial_k^h u(x):=\frac{u(x+he_k)-u(x)}h,\qquad |h|\ll1.
\end{equation*}
We recall that the action of the difference quotient on the product and the integral is similar to the derivative. Namely for $u,v:\R^n\to\R$ there holds
\begin{equation*}
\partial_k^h(uv)(x)=u(x+he_k)\partial_k^hv(x)+\partial_k^hu(x)v(x)
\end{equation*}
and
\begin{equation}\label{IPP}
\int u(x)\partial_k^hv(x)\,dx=-\int \partial_k^{-h}u(x)v(x)\,dx.
\end{equation}
Moreover, we will make use of the following well-known properties of the difference quotient, see e.g. \cite[Chp. 2 Lemma 4.6]{LU} or \cite[Lemma 2.7]{Mar}:
\begin{lem}\label{Lemma_Mar}
	Let $\Omega'$ be an open set compactly contained in $\Omega\subset\R^n$ and let $h_0:=\dist(\Omega',\R^n\setminus\Omega)$.
	\begin{enumerate}
		\item If $v\in W^{1,q}(\Omega)$ for some $q\geq1$, then for every $|h|\leq h_0$ there holds $\int_{\Omega'}|\phk v|^q\leq\int_\Omega|\partial_kv|^q$;
		\item If $v\in L^q(\Omega)$ for some $q>1$  and if there exists a constant $C>0$ such that $\|\phk v\|_{L^q(\Omega')}\leq C$ for every $|h|\leq h_0$, then $\partial_k v\in L^q(\Omega')$ and $\|\partial_kv\|_{L^q(\Omega')}\leq C$;
		\item If $v\in W^{1,q}(\Omega)$ for some $q>1$, then $\phk v\to\partial_kv$ strongly in $L^q(\Omega')$ for every $k\in\{1,\dots,n\}$.
	\end{enumerate}
	When $\Omega=\R^n$, points \textit{1.}-\textit{3.} hold also for $\Omega'=\R^n$ and any $h_0>0$.
\end{lem}

\subsection{Local estimates for the derivatives of \texorpdfstring{$\boldsymbol{\nabla\phi}$}{nabla phi}}\label{local_der}
Let $k\in\{1,\dots,n\}$ and 
\begin{equation}\label{zeta}
\zeta\in\begin{cases}
C^\infty_0(\R^n_\pm)&\quad k=1,\\
C^\infty_0(\R^n)&\quad k\geq2,
\end{cases}
\end{equation}
with $0<\zeta\leq 1$ in the interior of its support. In this section, we investigate at once \textit{local} estimates for the normal derivative of $\nabla\phi$, which corresponds to the case $k=1$, and for the tangential derivatives of $\nabla\phi$, i.e. the case $k\in\{2,\dots,n\}$. In both cases we aim to obtain the estimate
\begin{equation}\label{local_estimate_der}
\begin{split}
\int|\partial_k\nabla\phi|^2\zeta^2+\int|U_0+\nabla\phi|^{p-1}|\partial_k\nabla\phi|^2\zeta^2&\leq C\Big(\|U_0\|_{p+1}^{p+1}+\|U_0\|_{\alpha+1}^{\alpha+1}+\|\partial_kU_0\|_{p+1}^{p+1}\\
&\quad+\|\partial_kU_0\|_{\alpha+1}^{\alpha+1}+\|b\|_2^2+\iii b^2+\|\partial_kb\|_2^2\Big),
\end{split}
\end{equation}
where in the case $k=1$ the Lebesgue norms of $\partial_kU_0$ and $\partial_kb$ have to be understood on $\R^n_\pm$. To this purpose it is important to choose the test function $\zeta$ as in \eqref{zeta}, and hence, when the normal derivative is considered, $\supp\,\zeta$ is strictly contained in one of the two half-spaces. We stress the fact that the constant $C$ in \eqref{local_estimate_der} will depend on $\|\nabla\zeta\|_{\infty}$ in addition to the structural constants of the problem as detailed in Theorem \ref{Thm_der_phi}.
\vskip0.2truecm
We test \eqref{wf} with $\eta=\partial_k^{-h}\!\left(\phk\phi\,\zeta^2\right)$ for $h>0$ when $k\in\{2,\dots,n\}$ and $0<h<\dist(\supp\,\zeta,\Gamma)$ when $k=1$. Analogously to \eqref{b_Gauss} we apply Gau\ss{}'s theorem, obtaining
\begin{equation}\label{wf_der}
\int\epsilon_ff(U_0+\nabla\phi)\cdot\partial_k^{-h}\nabla\left(\phk\phi\,\zeta^2\right)+\int\epsilon_1\nabla\phi\cdot\partial_k^{-h}\nabla\left(\phk\phi\,\zeta^2\right)=\int b\,\partial_k^{-h}\!\left(\phk\phi\,\zeta^2\right).
\end{equation}
We begin by considering the linear terms. First by \eqref{IPP},
\begin{equation}\label{all_linear_terms}
\begin{split}
-\int\epsilon_1\nabla\phi\cdot\partial_k^{-h}\nabla\left(\phk\phi\,\zeta^2\right)&=\int\phk\left(\epsilon_1\nabla\phi\right)\ccdot\nabla\left(\phk\phi\,\zeta^2\right)\\
&=\int\epsilon_1(\cdot+he_k)|\nabla\phk\phi|^2\zeta^2+2\int\epsilon_1(\cdot+he_k)\nabla\phk\phi\cdot\phk\phi\,\zeta\nabla\zeta\\
&\quad+\int\phk\epsilon_1\nabla\phi\cdot\nabla\phk\phi\,\zeta^2+2\int\phk\epsilon_1\nabla\phi\cdot\phk\phi\,\zeta\nabla\zeta.
\end{split}
\end{equation}
We estimate the first term on the right in \eqref{all_linear_terms} from below by
\begin{equation}\label{linear_term_nablaphi2}
\int\epsilon_1(\cdot+he_k)|\nabla\phk\phi|^2\zeta^2\geq d\int|\nabla\phk\phi|^2\zeta^2,
\end{equation}
and we estimate from above all other terms. Here we use $\delta$-Young inequalities, Lemma \ref{Lemma_Mar}, and the estimate $|\phk\epsilon_1|\leq\|\partial_k\epsilon_1\|_{W^{1,\infty}(\R^n\setminus \Gamma)}$, which holds due to our choice of $h$:
\begin{equation}\label{linear_term_2}
\begin{split}
\left|\int\epsilon_1(\cdot+he_k)\nabla\phk\phi\cdot\phk\phi\,\zeta\nabla\zeta\right|&\leq\|\epsilon_1\|_\infty\|\nabla\zeta\|_\infty\left(\delta\int|\nabla\phk\phi|^2\zeta^2+C_\delta\|\phk\phi\,\zeta\|_2^2\right)\\
&\les\delta\int|\nabla\phk\phi|^2\zeta^2+C_\delta\|\nabla\phi\|_2^2,
\end{split}
\end{equation}
\begin{equation}\label{linear_term_3}
\left|\int\phk\epsilon_1\nabla\phi\cdot\nabla\phk\phi\,\zeta^2\right|\leq\|\epsilon_1\|_{W^{1,\infty}(\R^n\setminus \Gamma)}\left(\delta\int|\nabla\phk\phi|^2\zeta^2+C_\delta\|\nabla\phi\|_2^2\right),
\end{equation}
\begin{equation}\label{linear_term_4}
\left|\int\phk\epsilon_1\nabla\phi\cdot\phk\phi\,\zeta\nabla\zeta\right|\leq\|\epsilon_1\|_{W^{1,\infty}(\R^n\setminus \Gamma)}\|\nabla\zeta\|_\infty\left(\|\nabla\phi\|_2^2+\|\phk\phi\,\zeta\|_2^2\right)\les\|\nabla\phi\|_2^2.
\end{equation}
Notice that all terms on the right-hand sides of \eqref{linear_term_2}-\eqref{linear_term_4} either may be absorbed by \eqref{linear_term_nablaphi2} or depend just on $\|\nabla\phi\|_2$, which can be estimated by known quantities using Theorem \ref{Thm_phi}. We underline the fact that here, as well as in the following estimates, the constants hidden in the symbol $\les$ may also depend on $\|\nabla\zeta\|_{\infty}$.

The last term in \eqref{wf_der}, once integrated by parts as in \eqref{IPP}, can be easily estimated:
\begin{equation}\label{eq_L2_b_der}
\left|\int b\,\partial_k^{-h}\!\left(\phk\phi\zeta^2\right)\right|\leq\|\phk b\,\zeta\|_2\|\phk\phi\,\zeta\|_2\leq\|\partial_k b\|_2^2+\|\nabla\phi\|_2^2.
\end{equation}

The next step is to estimate the first term in \eqref{wf_der}. Following \cite{Mar}, we define
$$g:\R^n\times\R^n\to\R^n,\qquad g(x,\xi):=\epsilon_f(x)f(\xi)$$
and so for $u:\R^n\to\R^n$ one has
\begin{equation*}
\begin{split}
\phk g(x,u(x))&=\frac1h\int_0^1\frac d{dt}g(x+the_k,u(x)+th\phk u(x))\,dt\\
&=\int_0^1\!\left(\frac{\partial g}{\partial x_k}(x+the_k,u(x)+th\phk u(x))+\sum_{j=1}^n \frac{\partial g}{\partial \xi_j}(x+the_k,u(x)+th\phk u(x))\phk u_j(x)\!\right)\!dt\\
&=\int_0^1\!\partial_k\epsilon_f(x+the_k)f(u(x)+th\phk u(x))\,dt+\int_0^1\epsilon_f(x+the_k)J_f\!\left(u(x)+th\phk u(x)\right)\!\phk u(x)\,dt.
\end{split}
\end{equation*}
Hence, defining
$$Z_{th}:=(U_0+\nabla\phi)+th\phk(U_0+\nabla\phi)=(1-t)(U_0+\nabla\phi)+t(U_0+\nabla\phi)(\cdot+he_k),$$
we get for $u=U_0+\nabla\phi$
\begin{equation}\label{stime_der_NL}
\begin{split}
-\!\int\epsilon_ff(U_0+\nabla\phi)&\cdot\partial_k^{-h}\,\nabla\!\left(\phk\phi\,\zeta^2\right)=\int\!\!\left(\int_0^1\partial_k\epsilon_f(\cdot+the_k)f(\Zth)\,dt\right)\ccdot\left(\nabla\phk\phi\,\zeta^2+2\phk\phi\,\zeta\nabla\zeta\right)\\
&+\int\!\!\left(\int_0^1\epsilon_f(\cdot+the_k)J_f(\Zth)\phk(U_0+\nabla\phi)\,dt\right)\!\ccdot\left(\nabla\phk\phi\,\zeta^2+2\phk\phi\,\zeta\nabla\zeta\right)\\
&=:T_1+T_2+S_{11}+S_{12}+S_{21}+S_{22},
\end{split}
\end{equation}
using the discrete partial integration \eqref{IPP}. The $T$-terms (resp. $S$-terms) originate from the two (resp. four) products contained in the first (resp. second) term. Let us analyse all these terms separately, starting with $S_{21}$, since it contains the quantity which is going to enter on the left-hand side of the desired inequality \eqref{estimate_der_phi_tg}. Indeed, by means of (H0) and (F2.ii), one has
\begin{equation}\label{S21_big_term}
\begin{split}
S_{21}:=\int\!\!\left(\int_0^1\epsilon_f(\cdot+the_k)J_f(\Zth)\nabla\phk\phi\,dt\right)\ccdot\nabla\phk\phi\,\zeta^2\geq d\tlp\int\!\!\left(\int_0^1|\Zth|^{p-1}\,dt\right)|\nabla\phk\phi|^2\zeta^2.
\end{split}
\end{equation}
All the remaining terms will be estimated from above by means of (F2.i), either by the terms in the left-hand sides of \eqref{S21_big_term} and \eqref{linear_term_nablaphi2} multiplied by $\delta$, or by known quantities. For instance
\begin{equation*}
\begin{split}
|S_{11}|&\leq\|\epsilon_f\|_\infty\int\!\!\left(\tLp\int_0^1|\Zth|^{p-1}\,dt+\tLa\int_0^1|\Zth|^{\alpha-1}\,dt\right)|\phk U_0||\phk\nabla\phi|\zeta^2\\
&\les\delta\int\!\!\left(\int_0^1|\Zth|^{p-1}\,dt+\int_0^1|\Zth|^{\alpha-1}\right)|\nabla\phk\phi|^2\zeta^2\\
&\quad+C_\delta\int\!\!\left(\int_0^1|\Zth|^{p-1}\,dt+\int_0^1|\Zth|^{\alpha-1}\,dt\right)|\phk U_0|^2\zeta^2,\\
&\les\delta\int\!\!\left(\int_0^1|\Zth|^{p-1}\,dt+1\right)|\nabla\phk\phi|^2\zeta^2\\
&\quad+C_\delta\int\!\!\left(\int_0^1|\Zth|^{p-1}\,dt+\int_0^1|\Zth|^{\alpha-1}\,dt\right)|\phk U_0|^2\zeta^2,
\end{split}
\end{equation*}
having used the $\delta$-Young inequality in the first step and the estimates $x^{\alpha-1}\les x^{p-1}+1$ for $x\geq0$, since $\alpha\in(1,p]$, and $\|\zeta\|_\infty\leq1$ in the second step. In the first term we recognise the ``good'' quantities. Noticing that by the definition of $\Zth$ one has
\begin{equation}\label{Zth_estimate}
\int_0^1|\Zth|^{p-1}\,dt\les|U_0+\nabla\phi|^{p-1}+|U_0+\nabla\phi|^{p-1}(\cdot+he_k),
\end{equation}
and similarly for $\alpha$ replacing $p$, we get via H\"older and Young inequalities
\begin{equation*}
\begin{split}
\int\!\!\left(\int_0^1|\Zth|^{p-1}\,dt\right)|\phk U_0|^2\zeta^2&\les\|U_0+\nabla\phi\|_{p+1}^{p+1}+\|(U_0+\nabla\phi)(\cdot+he_k)\|_{p+1}^{p+1}+2\|\phk U_0\,\zeta\|_{p+1}^{p+1}\\
&\les2\|U_0+\nabla\phi\|_{p+1}^{p+1}+2\|\partial_kU_0\|_{p+1}^{p+1}\\
&\les\int|U_0+\nabla\phi|^{p-1}|\nabla\phi|^2+\|U_0\|_{p+1}^{p+1}+\|\partial_kU_0\|_{p+1}^{p+1},
\end{split}
\end{equation*}
having used \eqref{weighted_term_Lp+1_converse} in the last step. Similarly
\begin{equation*}
\begin{split}
\int\!\!\left(\int_0^1|\Zth|^{\alpha-1}\,dt\right)|\phk U_0|^2\zeta^2&\les\int|U_0+\nabla\phi|^{\alpha-1}|\phk U_0|^2\zeta^2+\int|U_0+\nabla\phi|^{\alpha-1}(\cdot+he_k)|\phk U_0|^2\zeta^2\\
&\les\|U_0\|_{\alpha+1}^{\alpha+1}+\|\phk U_0\,\zeta\|_{\alpha+1}^{\alpha+1}+\|\nabla\phi\|_{\alpha+1}^{\alpha+1}\\
&\les\|U_0\|_{\alpha+1}^{\alpha+1}+\|\partial_kU_0\|_{\alpha+1}^{\alpha+1}+\|\nabla\phi\|_{p+1}^{p+1}+\|\nabla\phi\|_2^2\\
&\les\|U_0\|_{\alpha+1}^{\alpha+1}+\|\partial_kU_0\|_{\alpha+1}^{\alpha+1}+\|U_0\|_{p+1}^{p+1}+\|\nabla\phi\|_2^2+\int|U_0+\nabla\phi|^{p-1}|\nabla\phi|^2.
\end{split}
\end{equation*}
Here we used $x^{\alpha+1}\les x^{p+1}+x^2$ for $x\geq0$ in the second inequality since $\alpha\in(1,p]$, and \eqref{nablaphi_p+1_pt2} in the last inequality to estimate the third term. All in all we get
\begin{equation*}
\begin{split}
|S_{11}|&\les\delta\int\!\!\left(\int_0^1|\Zth|^{p-1}\,dt+1\right)|\nabla\phk\phi|^2\zeta^2+C_\delta\Big(\|U_0\|_{\alpha+1}^{\alpha+1}+\|\partial_kU_0\|_{\alpha+1}^{\alpha+1}\\
&\quad+\|U_0\|_{p+1}^{p+1}+\|\partial_kU_0\|_{p+1}^{p+1}+\int|U_0+\nabla\phi|^{p-1}|\nabla\phi|^2+\|\nabla\phi\|_2^2\Big).
\end{split}
\end{equation*}
The terms $S_{12}$ and $S_{22}$ are estimated in the same way, obtaining
\begin{equation*}
|S_{12}|\les\|U_0\|_{\alpha+1}^{\alpha+1}+\|\partial_kU_0\|_{\alpha+1}^{\alpha+1}+\|U_0\|_{p+1}^{p+1}+\int|U_0+\nabla\phi|^{p-1}|\nabla\phi|^2+\|\nabla\phi\|_2^2
\end{equation*}
and
\begin{equation*}
\begin{split}
|S_{22}|&\les\delta\int\left(\int_0^1|\Zth|^{p-1}\,dt+1\right)|\nabla\phk\phi|^2\zeta^2+C_\delta\left(\|U_0\|_{p+1}^{p+1}+\int|U_0+\nabla\phi|^{p-1}|\nabla\phi|^2+\|\nabla\phi\|_2^2\right).
\end{split}
\end{equation*}
It remains to estimate the $T$-terms in \eqref{stime_der_NL}. Using (F1.i) we have
\begin{equation}\label{T_1}
\begin{split}
|T_1|&\leq\|\epsilon_f\|_{W^{1,\infty}(\R^n\setminus \Gamma)}\int\!\!\left(\Lp\int_0^1|\Zth|^p\,dt+\La\int_0^1|\Zth|^\alpha\,dt\right)|\nabla\phk\phi|\zeta^2.\\
\end{split}
\end{equation}
Concerning the first term, we exchange the integrals by Fubini in order to use H\"older and $\delta$-Young inequalities in the $x$-integral, restoring the original integration order by a second application of Fubini:
\begin{equation*}
\begin{split}
\int\!\!\bigg(\int_0^1|\Zth|^p\,&dt\bigg)|\nabla\phk\phi|\zeta^2=\int_0^1\left(\int|Z_{th}|^{\frac{p-1}2}|\nabla\phk\phi|\,\zeta\,\,|Z_{th}|^{\frac{p+1}2}\zeta\right)dt\\
&\leq\delta\int\!\!\left(\int_0^1|\Zth|^{p-1}\,dt\right)|\nabla\phk\phi|^2\zeta^2+C_\delta\int\!\!\left(\int_0^1|\Zth|^{p+1}\,dt\right)\\
&\les\delta\int\!\!\left(\int_0^1|\Zth|^{p-1}\,dt\right)|\nabla\phk\phi|^2\zeta^2+C_\delta\left(\|U_0\|_{p+1}^{p+1}+\int|U_0+\nabla\phi|^{p-1}|\nabla\phi|^2\right).
\end{split}
\end{equation*}
Note that in the last step we argued similarly to \eqref{Zth_estimate} with $p+1$ instead of $p-1$, and then applied \eqref{weighted_term_Lp+1_converse} to the resulting terms.

Analogously one may proceed with the second term in \eqref{T_1}:
\begin{equation*}
\begin{split}
\int\!\!\left(\int_0^1|\Zth|^\alpha\,dt\right)|\nabla\phk\phi|\zeta^2&\leq\delta\int\!\!\left(\int_0^1|\Zth|^{\alpha-1}\,dt\right)|\nabla\phk\phi|^2\zeta^2+C_\delta\int\!\!\left(\int_0^1|\Zth|^{\alpha+1}\,dt\right)\\
&\les\delta\int\!\!\left(\int_0^1|\Zth|^{p-1}\,dt+1\right)|\nabla\phk\phi|^2\zeta^2\\
&\quad+C_\delta\left(\|U_0\|_{\alpha+1}^{\alpha+1}+\|U_0\|_{p+1}^{p+1}+\|\nabla\phi\|_2^2+\int|U_0+\nabla\phi|^{p-1}|\nabla\phi|^2\right).
\end{split}
\end{equation*}
Hence, from \eqref{T_1} we find
\begin{equation*}
\begin{split}
|T_1|&\les\delta\int\!\!\left(\int_0^1|\Zth|^{p-1}\,dt+1\right)|\nabla\phk\phi|^2\zeta^2\\
&\quad+C_\delta\left(\|U_0\|_{\alpha+1}^{\alpha+1}+\|U_0\|_{p+1}^{p+1}+\|\nabla\phi\|_2^2+\int|U_0+\nabla\phi|^{p-1}|\nabla\phi|^2\right).
\end{split}
\end{equation*}
Finally, a similar bound independent of $\delta$ can be established for $|T_2|$ since it does not involve $\nabla\phk\phi\,$:
\begin{equation*}
|T_2|\les\|U_0\|_{\alpha+1}^{\alpha+1}+\|U_0\|_{p+1}^{p+1}+\|\nabla\phi\|_2^2+\int|U_0+\nabla\phi|^{p-1}|\nabla\phi|^2.
\end{equation*}
Therefore, gathering the estimates we obtained for the terms in \eqref{all_linear_terms}, \eqref{eq_L2_b_der}, and \eqref{stime_der_NL}, we infer
\begin{equation}\label{eq_quasi_finale}
\begin{split}
(1-C_1\delta)&\int|\nabla\phk\phi|^2\zeta^2+(1-C_2\delta)\int\!\!\left(\int_0^1|\Zth|^{p-1}\,dt\right)|\nabla\phk\phi|^2\zeta^2\\
&\les C_\delta\Big(\int|U_0+\nabla\phi|^{p-1}|\nabla\phi|^2+\|\nabla\phi\|_2^2+\|U_0\|_{\alpha+1}^{\alpha+1}\\
&\quad+\|\partial_kU_0\|_{\alpha+1}^{\alpha+1}+\|U_0\|_{p+1}^{p+1}+\|\partial_kU_0\|_{p+1}^{p+1}+\|b\|_2^2+\iii b^2+\|\partial_kb\|_2^2\Big)
\end{split}
\end{equation}
for suitable constants $C_1,C_2>0$. Choosing a sufficiently small $\delta$, and since the second term on the left in \eqref{eq_quasi_finale} is positive, one deduces that $\int|\nabla\phk\phi|^2\zeta^2$ is uniformly bounded with respect to $h$ and that the constant depends only on the $W^{1,\infty}$-norm of $\zeta$. Hence, by Lemma \ref{Lemma_Mar} one has $\nabla\partial_k\phi\in L^2(\supp\,\zeta)$ (with the $L^2$-norm bounded by the same right-hand side as in \eqref{eq_quasi_finale}) and $\nabla\phk\phi\to\nabla\partial_k\phi$ in $L^2(\supp\,\zeta)$ as $h\to0$, so pointwise a.e. in $\supp\,\zeta$. Moreover, thanks to the continuity with respect to translations in the $L^{p+1}$-norm, one has
$$Z_{th}=(1-t)(U_0+\nabla\phi)+t(U_0+\nabla\phi)(\cdot+he_k)\to U_0+\nabla\phi\quad\mbox{a.e. in}\;\R^n\quad\!\mbox{as}\;\,h\to0.$$
Hence, applying Fubini's theorem and Fatou's Lemma, one gets
\begin{equation*}
\begin{split}
\lim_{h\to0}\int\left(\int_0^1|Z_{th}|^{p-1}\,dt\right)|\nabla\phk\phi|^2\zeta^2&\geq\int_0^1\left(\int\lim_{h\to0}|Z_{th}|^{p-1}|\nabla\phk\phi|^2\zeta^2\right)\,dt\\
&=\int|U_0+\nabla\phi|^{p-1}|\nabla\partial_k\phi|^2\zeta^2,
\end{split}
\end{equation*}
and therefore from \eqref{eq_quasi_finale}, once estimate \eqref{estimate_phi} found in Theorem \ref{Thm_phi} is applied, one infers \eqref{local_estimate_der}.

\subsection{Tangential derivatives \texorpdfstring{$\boldsymbol{\partial_k \nabla \phi, k\geq2}$}{k>=2}: estimate \eqref{estimate_der_phi_tg}}\label{Sec_tan_der}
Let $\zeta_0\in C^\infty_0(\R^n)$ so that $\zeta_0\equiv1$ on $B_1(0)$ and $\supp\,\zeta_0\subset B_2(0)$, and let $(\zeta_j)_j$ be the sequence of test functions defined as $\zeta_j:=\zeta_0\big(\tfrac\cdot j\big)$. Then the supports of $(\zeta_j)_j$ grow to cover the whole $\R^n$, $\zeta_j(x)\to 1$ for all $x\in\R^n$ as $j\to+\infty$, and $\|\nabla\zeta_j\|_{W^{1,\infty}}$ are uniformly bounded. Hence an application of Fatou's Lemma to \eqref{local_estimate_der} with $\zeta=\zeta_j$ yields \eqref{estimate_der_phi_tg}.

Note that for the estimate of tangential derivatives of $\nabla\phi$ introducing the factor $\zeta$ in the test function $\eta$ is actually not necessary. Indeed, testing with $\eta=\partial_k^{-h}\left(\phk\phi\right)$ is allowed since Lemma \ref{Lemma_Mar} holds also for $\Omega=\Omega'=\R^n$. However, for the sake of a concise presentation in Section \ref{Sec_norm_der} we use a test function which can be used for all derivatives.

\subsection{Normal derivative \texorpdfstring{$\boldsymbol{\partial_{11}\phi}$}{partial11phi}: estimate \eqref{estimate_der_phi_n}}\label{Sec_norm_der}
As mentioned in Remark \ref{Rmk_norm_der}, we are only left with the estimate for the first component of the normal derivative $\partial_1\nabla\phi$. Indeed, by the argument detailed in Section \ref{local_der}, for a $\zeta\in C^\infty_0(\R^n_\pm)$ one infers the existence of $\partial_1\nabla\phi$ in $L^2(\supp\,\zeta)$ together with the estimate
\begin{equation}\label{pa1nabla}
\int|\partial_1\nabla\phi|^2\zeta^2+\int|U_0+\nabla\phi|^{p-1}|\partial_1\nabla\phi|^2\zeta^2\leq K_1,
\end{equation}
where the constant $K_1$ depends on the $L^{p+1}$- and $L^{\alpha+1}$-norms of $U_0,\,\partial_1U_0,\,\partial_1b$ and on $\iii b$ in the half-space which contains $\supp\,\zeta$. By the arbitrariness of $\zeta\in C^\infty_0(\R^n_\pm)$, one deduces that $\partial_1\nabla\phi\in L^2_{loc}(\R^n_\pm)^n$. Since $\nabla \phi\in H^1_{loc}(\R^n_\pm)^n$, it is easy to verify that $\partial_{1k}\phi=\partial_{k1}\phi$ for all $k=1,\dots,n$ as functions in $L^2(\R^n_\pm)$ and so also a.e. in $\R^n_\pm$. Hence, decomposing the domain of integration as $\R^n=\R^n_+\cup\R^n_-$, we get for $k=2,\dots,n$
\begin{equation}\label{E:Dk1}
\begin{split}
&\int|\partial_{1k}\phi|^2+\int|U_0+\nabla\phi|^{p-1}|\partial_{1k}\phi|^2=\int|\partial_{k1}\phi|^2+\int|U_0+\nabla\phi|^{p-1}|\partial_{k1}\phi|^2\leq\widetilde C_k,
\end{split}
\end{equation}
where $\widetilde C_k$ is the right-hand side of \eqref{estimate_der_phi_tg}. Hence we conclude that $\partial_{1k}\phi=\partial_{k1}\phi$ in $L^2(\R^n)$, which means that $\partial_{1k}\phi$, $k=2,\dots,n$ are estimated in Section \ref{Sec_tan_der}.

The trick in \eqref{E:Dk1} is of course not applicable for $\partial_{11}\phi$. Moreover, note that in the estimate we get from \eqref{pa1nabla}, i.e.
\begin{equation}\label{estimate_der_phi_n_loc}
\int|\partial_{11}\phi|^2\zeta^2+\int|U_0+\nabla\phi|^{p-1}|\partial_{11}\phi|^2\zeta^2\leq K_1,
\end{equation}
the constant $K_1$ \textit{does depend on} $\|\nabla\zeta\|_{\infty}$. Hence, when one defines a sequence of uniformly bounded test functions $(\zeta_j)_j$ similarly to Section \ref{Sec_tan_der}, namely such as their support grows to cover e.g. the subspace $\R^n_+$, and $\zeta_j\to1$ pointwise on $\R^n_+$, then $\|\nabla\zeta_j\|_{\infty}\to \infty$ as $j\to \infty$ because $\dist(\supp\,\zeta_j,\Gamma)\to 0$. This implies that simply applying Fatou's Lemma to \eqref{estimate_der_phi_n_loc}, as we did for the tangential derivatives, is insufficient. The aim of this section is thus to obtain an estimate for $\partial_{11}\phi$ similar to \eqref{estimate_der_phi_n_loc} but independent of $\|\nabla\zeta\|_{\infty}$.
\vskip0.2truecm
First, we claim that equation \eqref{eq} is actually satisfied \textit{pointwise a.e.} in $\R^n_\pm$. Indeed, from the weak formulation \eqref{wf} with $\eta\in C^\infty_0(\R^n_\pm)$ and from \eqref{b_Gauss} with such $\eta$ instead of $\phi$ we get by partial integration
\begin{equation}\label{wf_Gauss}
-\int_{\R^n_\pm}\nabla\cdot\left(\epsilon_ff(U_0+\nabla\phi)+\epsilon_1\nabla\phi\right)\eta=\int_{\R^n_\pm} b\eta,
\end{equation}
where the term on the left-hand side is well-defined. Indeed, concerning the second summand, the functions $\partial_k\left(\epsilon_1\partial_k\phi\right)\in L^1_{loc}(\R^n_\pm)$ because of \eqref{estimate_der_phi_n_loc} for $k=1$ and \eqref{estimate_der_phi_tg} for $k\in\{2,\dots,n\}$. For the first summand we have\footnote{We denote by $J_f^{(k,j)}$ the $(k,j)$-element of the matrix $J_f$ and by $J_f^{(k,\cdot)}$ its $k$-th row.},
\begin{equation*}
\partial_k\left(\epsilon_ff_k(U_0+\nabla\phi)\right)=\partial_k\epsilon_ff_k(U_0+\nabla\phi)+\epsilon_fJ_f^{(k,\cdot)}(U_0+\nabla\phi)\cdot\partial_k(U_0+\nabla\phi),
\end{equation*}
so by (F1.i)
\begin{equation*}
\begin{split}
\left|\int_{\supp\,\eta}\!\partial_k\epsilon_ff_k(U_0+\nabla\phi)\right|&\les\|\epsilon_f\|_{W^{1,\infty}(\R^n\setminus \Gamma)}\left(\int_{\supp\,\eta}|U_0+\nabla\phi|^p+\int_{\supp\,\eta}|U_0+\nabla\phi|^\alpha\right)\\
&\les\int_{\supp\,\eta}|U_0+\nabla\phi|^p+|\supp\,\eta|\les\|U_0+\nabla\phi\|_{p+1}^{p+1}+|\supp\,\eta|,
\end{split}
\end{equation*}
while by (F2.i)
\begin{equation*}
\begin{split}
\Big|\int_{\supp\,\eta}\!\epsilon_fJ_f^{(k,\cdot)}(U_0+\nabla\phi)\cdot&\partial_k(U_0+\nabla\phi)\Big|\les\|\epsilon_f\|_{\infty}\int_{\supp\,\eta}\left(|U_0+\nabla\phi|^{p-1}+1\right)\left(|\partial_kU_0|+|\partial_k\nabla\phi|\right)\\
&\les\int|U_0+\nabla\phi|^{p+1}+\int_{\supp\,\eta}|U_0+\nabla\phi|^{p-1}|\partial_k\nabla\phi|^2+\int_{\supp\,\eta}|\partial_k\nabla\phi|^2\\
&\quad+\int_{\supp\,\eta}|\partial_kU_0|^{p+1}+\int_{\supp\,\eta}|\partial_kU_0|^2+|\supp\,\eta|.
\end{split}
\end{equation*}
Here we used the trivial inequality $t^{p-1}\les t^p+1$ for $t\geq0$, then wrote $p=\frac{p+1}2+\frac{p-1}2$ and applied the Cauchy-Schwarz inequality.

Then, by the arbitrariness of $\eta\in C^\infty_0(\R^n_\pm)$ and invoking the fundamental lemma of calculus of variations, the fact that the solution $\nabla\phi$ satisfies equation \eqref{eq} a.e. follows from \eqref{wf_Gauss}.
\vskip0.2truecm

Next, we prove the interface condition \eqref{IFCs_eq} in the trace sense of Definition \ref{def_trace}. As $\nabla\cdot \cD(U_0+\nabla \phi)=0$ a.e. in $\R^n$, we get
$$(T_\pm\cD)[\psi]=\mp \int_{\R^n_\pm}\cD\cdot \nabla \hat{\psi},$$
where $\hat{\psi}\in C^\infty_0(\R^n)$ is such that $\hat{\psi}|_\Gamma=\psi$. Hence we have
$$(T_+\cD-T_-\cD)[\psi]=-\int_{\R^n}\cD\cdot \nabla \hat{\psi}=0$$
using the definition of the weak solution and the fact that $C^\infty_0(\R^n)\subset \fD_{2,p+1}$.

\vskip0.2truecm
In the final part of the  proof we estimate the normal derivative $\partial_{11}\phi$. Let us now choose a half-space, say $\R^n_+$, and multiply equation \eqref{eq} by $\partial_{11}\phi\,\zeta$, with $\zeta\in C^\infty_0(\R^n_+)$. Integrating over $\R^n_+$, we get
\begin{equation}\label{eq_pw_der11}
\int_{\R^n_+}\nabla\cdot\left(\epsilon_ff(U_0+\nabla\phi)\right)\partial_{11}\phi\,\zeta+\int_{\R^n_+}\nabla\cdot\left(\epsilon_1\nabla\phi\right)\partial_{11}\phi\,\zeta=-\int_{\R^n_+}b\,\partial_{11}\phi\,\zeta.
\end{equation}
We have
\begin{equation}\label{eq_pw_der11_LIN}
\int_{\R^n_+}\nabla\cdot\left(\epsilon_1\nabla\phi\right)\partial_{11}\phi\,\zeta=\int_{\R^n_+}\epsilon_1|\partial_{11}\phi|^2\zeta+\sum_{k=2}^n\int_{\R^n_+}\epsilon_1\partial_{kk}\phi\,\partial_{11}\phi\,\zeta+\sum_{k=1}^n\int_{\R^n_+}\partial_k\epsilon_1\partial_k\phi\,\partial_{11}\phi\,\zeta
\end{equation}
and for $k\in\{2,\dots,n\}$
\begin{equation*}
\left|\int_{\R^n_+}\epsilon_1\partial_{kk}\phi\,\partial_{11}\phi\,\zeta\right|\les\delta\intpiu|\partial_{11}\phi|^2\zeta+C_\delta\intpiu|\partial_{kk}\phi|^2,
\end{equation*}
while for $k\in\{1,\dots,n\}$
\begin{equation*}
\left|\int_{\R^n_+}\partial_k\epsilon_1\partial_k\phi\,\partial_{11}\phi\,\zeta\right|\les\delta\intpiu|\partial_{11}\phi|^2\zeta+C_\delta\intpiu|\nabla\phi|^2.
\end{equation*}
We can decompose the first term in \eqref{eq_pw_der11} as
\begin{equation}\label{eq_pw_der11_NONLIN}
\begin{split}
\int_{\R^n_+}&\nabla\cdot\left(\epsilon_ff(U_0+\nabla\phi)\right)\partial_{11}\phi\,\zeta=\sum_{k=1}^n\intpiu\epsilon_fJ_f^{(k,\cdot)}(U_0+\nabla\phi)\,\partial_k\nabla\phi\,\partial_{11}\phi\,\zeta\\
&+\sum_{k=1}^n\intpiu\epsilon_fJ_f^{(k,\cdot)}(U_0+\nabla\phi)\,\partial_kU_0\,\partial_{11}\phi\,\zeta+\sum_{k=1}^n\intpiu\partial_k\epsilon_ff_k(U_0+\nabla\phi)\,\partial_{11}\phi\,\zeta.
\end{split}
\end{equation}
The last two terms in \eqref{eq_pw_der11_NONLIN} are easy to handle. For $k\in\{1,\dots,n\}$ by (F2.i)
\begin{equation*}
\begin{split}
\left|\intpiu\epsilon_fJ_f^{(k,\cdot)}(U_0+\nabla\phi)\,\partial_kU_0\,\partial_{11}\phi\,\zeta\right|&\les\intpiu|U_0+\nabla\phi|^{p-1}|\partial_kU_0||\partial_{11}\phi|\zeta\\
&\quad+\intpiu|U_0+\nabla\phi|^{\alpha-1}|\partial_kU_0||\partial_{11}\phi|\zeta\\
&:=T_p+T_\alpha
\end{split}
\end{equation*}
where
\begin{equation*}
\begin{split}
T_p&\leq C_\delta\intpiu|U_0+\nabla\phi|^{p-1}|\partial_kU_0|^2+\delta\intpiu|U_0+\nabla\phi|^{p-1}|\partial_{11}\phi|^2\zeta^2\\
&\les C_\delta\left(\int|U_0+\nabla\phi|^{p+1}+\intpiu|\partial_kU_0|^{p+1}\right)+\delta\intpiu|U_0+\nabla\phi|^{p-1}|\partial_{11}\phi|^2\zeta^2,
\end{split}
\end{equation*}
and
\begin{equation*}
\begin{split}
T_\alpha&\leq C_\delta\left(\int|U_0+\nabla\phi|^{\alpha+1}+\intpiu|\partial_kU_0|^{\alpha+1}\right)+\delta\intpiu|U_0+\nabla\phi|^{\alpha-1}|\partial_{11}\phi|^2\zeta^2\\
&\les C_\delta\left(\|U_0\|_{\alpha+1}^{\alpha+1}+\|\nabla\phi\|_{p+1}^{p+1}+\|\nabla\phi\|_2^2+\intpiu|\partial_kU_0|^{\alpha+1}\right)\\
&\quad+\delta\left(\intpiu|U_0+\nabla\phi|^{p-1}|\partial_{11}\phi|^2\zeta^2+\intpiu|\partial_{11}\phi|^2\zeta^2\right)\\
&\les C_\delta\left(\|U_0\|_{\alpha+1}^{\alpha+1}+\|U_0\|_{p+1}^{p+1}+\int|U_0+\nabla\phi|^{p-1}|\nabla\phi|^2+\|\nabla\phi\|_2^2+\intpiu|\partial_kU_0|^{\alpha+1}\right)\\
&\quad+\delta\left(\intpiu|U_0+\nabla\phi|^{p-1}|\partial_{11}\phi|^2\zeta^2+\intpiu|\partial_{11}\phi|^2\zeta^2\right),
\end{split}
\end{equation*}
where in the last step we applied \eqref{nablaphi_p+1_pt2}. Hence, using \eqref{weighted_term_Lp+1_converse} we can estimate
\begin{equation*}
\begin{split}
\bigg|\intpiu\epsilon_f&J_f^{(k,\cdot)}(U_0+\nabla\phi)\cdot\partial_kU_0\,\partial_{11}\phi\,\zeta\bigg|\les C_\delta\bigg(\|U_0\|_{p+1}^{p+1}+\|U_0\|_{\alpha+1}^{\alpha+1}+\int|U_0+\nabla\phi|^{p-1}|\nabla\phi|^2+\|\nabla\phi\|_2^2\\
&\quad+\intpiu|\partial_kU_0|^{p+1}+\intpiu|\partial_kU_0|^{\alpha+1}\bigg)+\delta\left(\intpiu|U_0+\nabla\phi|^{p-1}|\partial_{11}\phi|^2\zeta^2+\intpiu|\partial_{11}\phi|^2\zeta^2\right).
\end{split}
\end{equation*}
Similarly, by (F1.i) one has
\begin{equation*}
\begin{split}
\left|\intpiu\partial_k\epsilon_ff_k(U_0+\nabla\phi)\,\partial_{11}\phi\,\zeta\right|&\les\intpiu|U_0+\nabla\phi|^p|\partial_{11}\phi|\zeta+\intpiu|U_0+\nabla\phi|^\alpha|\partial_{11}\phi|\zeta\\
&\leq C_\delta\int|U_0+\nabla\phi|^{p+1}+\delta\intpiu|U_0+\nabla\phi|^{p-1}|\partial_{11}\phi|^2\zeta\\
&\quad+C_\delta\int|U_0+\nabla\phi|^{\alpha+1}+\delta\intpiu|U_0+\nabla\phi|^{\alpha-1}|\partial_{11}\phi|^2\zeta\\
&\les\delta\Big(\intpiu|U_0+\nabla\phi|^{p-1}|\partial_{11}\phi|^2\zeta+\intpiu|\partial_{11}\phi|^2\zeta\Big)\\
&\quad+C_\delta\Big(\int|U_0+\nabla\phi|^{p+1}+\int|\nabla\phi|^2+\|U_0\|_{p+1}^{p+1}+\|U_0\|_{\alpha+1}^{\alpha+1}\Big).
\end{split}
\end{equation*}
The first term on the right hand side in \eqref{eq_pw_der11_NONLIN} can be rewritten as follows:
\begin{equation}\label{Jkj_1}
\begin{split}
\sum_{k=1}^n&\intpiu\epsilon_fJ_f^{(k,\cdot)}(U_0+\nabla\phi)\,\partial_k\nabla\phi\,\partial_{11}\phi\,\zeta\\
&=\intpiu\epsilon_fJ_f^{(1,1)}(U_0+\nabla\phi)|\partial_{11}\phi|^2\zeta+\sum_{\substack{k,j=1\\(k,j)\neq(1,1)}}^n\intpiu \epsilon_fJ_f^{(k,j)}(U_0+\nabla\phi)\partial_{kj}\phi\,\partial_{11}\phi\,\zeta\\
&=\sum_{k,j=1}^n\intpiu\epsilon_fJ_f^{(k,j)}(U_0+\nabla\phi)\partial_{1k}\phi\,\partial_{1j}\phi\,\zeta\\
&\quad+\sum_{\substack{k,j=1\\(k,j)\neq(1,1)}}^n\intpiu\epsilon_fJ_f^{(k,j)}(U_0+\nabla\phi)\left(\partial_{kj}\phi\,\partial_{11}\phi-\partial_{1k}\phi\,\partial_{1j}\phi\right)\zeta.
 \end{split}
\end{equation}
Since
\begin{equation*}
\begin{split}
\sum_{k,j=1}^n\intpiu\epsilon_fJ_f^{(k,j)}(U_0+\nabla\phi)\,\partial_{1k}\phi\,\partial_{1j}\phi\,\zeta&=\intpiu\epsilon_f\left(J_f(U_0+\nabla\phi)\,\partial_1\nabla\phi\right)\ccdot\partial_1\nabla\phi\,\zeta\\
&\geq d\tlp\intpiu|U_0+\nabla\phi|^{p-1}|\partial_1\nabla\phi|^2\zeta\\
&\geq d\tlp\intpiu|U_0+\nabla\phi|^{p-1}|\partial_{11}\phi|^2\zeta,
\end{split}
\end{equation*}
by (F2.ii), we may recognise the term we want to keep in the final estimate. The second term on the right in \eqref{Jkj_1} is estimated from above as usual by (F2.i) and the $\delta$-Young inequality:
\begin{equation*}
\begin{split}
\left|\intpiu\epsilon_fJ_f^{(k,j)}(U_0+\nabla\phi)\,\partial_{kj}\phi\,\partial_{11}\phi\,\zeta\right|&\les \delta\intpiu|U_0+\nabla\phi|^{p-1}|\partial_{11}\phi|^2\zeta+\delta\intpiu|\partial_{11}\phi|^2\zeta\\
&\quad+C_\delta\int|U_0+\nabla\phi|^{p-1}|\partial_{kj}\phi|^2+C_\delta\int|\partial_{kj}\phi|^2\zeta,
\end{split}
\end{equation*}
\begin{equation*}
\begin{split}
\left|\intpiu\epsilon_fJ_f^{(k,j)}(U_0+\nabla\phi)\,\partial_{1k}\phi\,\partial_{1j}\phi\,\zeta\right|&\les \int|U_0+\nabla\phi|^{p-1}|\partial_{1k}\phi|^2+\int|\partial_{1k}\phi|^2.
\end{split}
\end{equation*}
This concludes the estimates for all terms in \eqref{eq_pw_der11_NONLIN}. Finally, the last term in \eqref{eq_pw_der11} is straightforward:
\begin{equation*}
\left|-\intpiu b\,\partial_{11}\phi\,\zeta\right|\leq\delta\|\partial_{11}\phi\|_2^2+C_\delta\|b\|_2^2.
\end{equation*}

Recalling inequality \eqref{estimate_der_phi_tg}, from equation \eqref{eq_pw_der11} we may thus conclude (by choosing a sufficiently small $\delta$)
\begin{equation}\label{local_norm_indep_grad}
\begin{split}
\intpiu|\partial_{11}\phi|^2\zeta+\intpiu|U_0+\nabla\phi|^{p-1}|\partial_{11}\phi|^2\zeta\leq C_1^+,
\end{split}
\end{equation}
where the constant $C_1^+$ has the same dependencies as in Theorem \ref{Thm_der_phi}. Notice that now the constant $C_1^+$ does not depend on $\|\zeta\|_{W^{1,\infty}(\R^n_+)}$ anymore, but only on $\|\zeta\|_\infty$. Hence, similarly to Section \ref{Sec_tan_der}, one may consider a sequence of test functions $(\zeta_j)_j\subset C^\infty_0(\R^n_+)$ with a uniformly bounded $L^\infty$-norm, which pointwise converge to $1$ and the supports of which grow to cover the subspace $\R^n_+$. Then, applying Fatou's Lemma to \eqref{local_norm_indep_grad} with $\zeta=\zeta_j$, we finally get the desired inequality \eqref{estimate_der_phi_n}. The argument for the halfspace $\R^n_-$ is completely analogous.

\section{Application: enforcing \texorpdfstring{$\boldsymbol{\nabla\cdot\cD(\cdot,\cE)=0}$}{divD(E)=0} at \texorpdfstring{$\boldsymbol{t=0}$}{t=0} in the Maxwell problem for \texorpdfstring{$\boldsymbol{n=2}$}{n=2}}\label{Section_application}

Let us return back to the Maxwell setting described in the introduction and apply the theory to problem \eqref{E0}-\eqref{initial_div}-\eqref{Max_IFCs_korr}. We recall that we are considering two different materials of dielectric kind which are separated by the interface $\Gamma=\{x_1=0\}$. Within the whole section we restrict ourselves to the case $n=2$ and assume that the materials are independent of $x_2,\,x_3$. In addition, we choose the special case of the Kerr nonlinear dependence $\cD=\cD(\cE)$ as in \eqref{Kerr}. With the notation of \eqref{relationD(E)}, one has
\begin{equation*}
f(v)=|v|^2v,\qquad\epsilon_1(x)=\epsilon_1(x_1),\qquad\epsilon_f(x)=\epsilon_3(x_1).
\end{equation*}
Of course, we require that the material coefficients $\epsilon_1$, $\epsilon_3$ fulfil assumptions (H0)-(H1).

\subsection{Linear problem and numerical computation of the eigenvalues}\label{Section_linear}
We start by analysing the \textit{linear} Maxwell problem, i.e. when $\cD=\cD_{lin}:=\epsilon_1(x)\cE$ or, in other words, the case $\epsilon_3\equiv0$. In this setting one may look for solutions of the form
\begin{equation}\label{modes}
\cE(x,t)=\begin{pmatrix}\varphi_1(x_1)\\ \varphi_2(x_1)\\ \varphi_3(x_1)\end{pmatrix}\e^{\ri(kx_2-\omega t)}+\mbox{c.c.},
\qquad\cH(x,t)=\begin{pmatrix}\psi_1(x_1)\\ \psi_2(x_1)\\ \psi_3(x_1)\end{pmatrix}\e^{\ri(kx_2-\omega t)}+\mbox{c.c.}
\end{equation}
with $k\in\R$, $\omega\in\R\setminus\{0\}$ and $\varphi,\psi:\R\to\C^3$. We study \textit{transverse magnetic} (TM) modes, that is $\varphi=\left(\varphi_1,\varphi_2,0\right)^\trans$ and $\psi=\left(0,0,\psi_3\right)^\trans$. Note that in the special case of piecewise constant $\epsilon_1$, i.e. $\epsilon_1=\epsilon_1^+\chi_{\R_+}+\epsilon_1^-\chi_{\R_-}$ with $\epsilon_1^\pm\in \R$, one can show that all solutions of the form \eqref{modes} with $\varphi,\psi\in L^2(\R)^3$ are TM modes, see \cite{BDPW}. 

Arranging the nonvanishing components in the vector $w(x_1):=\left(\varphi_1(x_1),\varphi_2(x_1),\psi_3(x_1)\right)^\trans$, one sees that, in order for $(\cE,\cH)$ to be a solution of the linear Maxwell problem, $w$ must satisfy the eigenvalue problem
\begin{equation}\label{L_selfadj_eq}
\begin{cases}
L(k,\omega)w=0 & \mbox{in}\quad\R\setminus\{0\},\\
\left\llbracket\epsilon_1w_1\right\rrbracket=\left\llbracket w_2\right\rrbracket=\left\llbracket w_3\right\rrbracket=0,
\end{cases}
\end{equation}
where the operator $L(k,\omega):D(L)\to L^2(\R)^3$ is given by
\begin{equation}\label{L_k}
L(k,\omega)w:=\begin{pmatrix}\epsilon_1\omega w_1+kw_3\\\epsilon_1\omega w_2+\ri w'_3\\kw_1+\ri w_2'+\omega w_3\end{pmatrix}
\end{equation}
and its domain
\begin{equation*}
D(L) := \big\{w\in L^2(\R)^3\,\big|\,w_2,w_3 \in H^1(\R)\big\}
\end{equation*}
is dense in $L^2(\R)^3$. Note that, because $\omega\neq0$, an element $w\in D(L)$ for which $L(k,\omega)w=0$ also satisfies the interface conditions, since $w_2,w_3\in H^1(\R)\hookrightarrow C(\R)$ and $\epsilon_1\omega w_1+kw_3=0$ implies that $\epsilon_1w_1$ is also continuous. Moreover, the first two equations of $L(k,\omega)w = 0$ imply that
\begin{equation}\label{divDlin=0}
\left(\epsilon_1 w_1\right)' + \epsilon_1 \ri kw_2 =0
\end{equation}
in $\R \setminus \{0\}$, which means that $\nabla\cdot\cD_{lin}=0$ holds.

Wave numbers $k$ for which $w\in D(L)$ exists such that \eqref{L_selfadj_eq} holds are called \textit{eigenvalues}. The dependence $k=k(\omega)$ is called the \textit{dispersion relation}. In general there can be more than one eigenvalue for a given $\omega$. Note that for the pointwise constant case the eigenvalue is unique (up to the sign), see below.
\begin{remark}\label{real_mag_real}
One can show that it is always possible to find a solution $w$ such that $w_1$ and $w_3$ are real and $w_2$ is imaginary. Indeed, if one substitutes the ansatz $(\widetilde w_1,\ri\widetilde w_2, \widetilde w_3)$ in \eqref{L_selfadj_eq}, one obtains a real problem for $(\widetilde w_1,\widetilde w_2, \widetilde w_3)$.
\end{remark}

One may also rewrite problem \eqref{L_selfadj_eq}-\eqref{L_k} as a second-order ODE on $\R^+$ and $\R^-$. Indeed,
\begin{align*}
w_3'' &= \ri\epsilon_1'\omega w_2+\ri\epsilon_1\omega w_2' \\
&= \frac{\epsilon_1'}{\epsilon_1} w_3'-\epsilon_1\omega\left(\omega w_3 + k w_1\right) \\
&= \frac{\epsilon_1'}{\epsilon_1}w_3'-\epsilon_1\omega^2 w_3 + k^2 w_3.
\end{align*}  
From $\llbracket w_2 \rrbracket= 0$ we deduce $\left\llbracket\frac{w_3'}{\epsilon_1}\right\rrbracket= 0$. Therefore, one needs to solve the eigenvalue problem
\begin{equation}\label{eigenvalue_problem_numerics}
\begin{cases}
-w_3''+\frac{\epsilon_1'}{\epsilon_1}w_3'-\epsilon_1\omega^2 w_3=- k^2 w_3&\mbox{in}\quad\R\setminus\{0\},\\
\llbracket w_3 \rrbracket=\left\llbracket\tfrac{w_3'}{\epsilon_1 }\right\rrbracket=0
\end{cases}
\end{equation}
and, once $w_3$ is obtained, one may get the remaining components of $w$ by the relations $w_1 = -\frac k{\omega\epsilon_1} w_3$ and $w_2 = - \frac{\ri}{\omega\epsilon_1} w'_3$. We also see that the interface conditions $\llbracket \epsilon_1 w_1 \rrbracket = \llbracket w_2 \rrbracket = 0$ at $x_1=0$ are satisfied if $w_3$ solves \eqref{eigenvalue_problem_numerics}.

When $\epsilon_1=\epsilon_1^-\chi_{\R_-}+\epsilon_1^+\chi_{\R_+}$ with constants $\epsilon_1^\pm\in\R$, equation \eqref{eigenvalue_problem_numerics} reduces to a second-order ODE with constant coefficients and one handily infers the existence of exponentially decaying solutions of the differential equation and the first interface condition in \eqref{eigenvalue_problem_numerics} of the form
\begin{equation*}
w_3(x_1) = \begin{cases}
C\,\e^{-\sqrt{\lambda_{\omega,k}^-}x_1} & x_1\leq0,\\
C\,\e^{\sqrt{\lambda_{\omega,k}^+}x_1} & x_1 >0,
\end{cases}
\end{equation*}
where $C\in\R\setminus\{0\}$ if $\lambda_{\omega,k}^\pm:=k^2-\epsilon_1^\pm\omega^2>0$. In order for the second interface condition to be satisfied too, one needs
\begin{equation}\label{disp_rel}
-\frac{\epsilon_1^+}{\epsilon_1^-}=\sqrt{\frac{\lambda_{\omega,k}^+}{\lambda_{\omega,k}^-}},
\end{equation}
which in turn yields the dispersion relation 
\begin{equation*}
k^2=\omega^2\frac{\epsilon_1^+\epsilon_1^-}{\epsilon_1^++\epsilon_1^-}.
\end{equation*}
However, \eqref{disp_rel} implies that $\epsilon_1^+$ and $\epsilon_1^-$ should have opposite signs, which is not allowed by assumption (H0). 

Therefore, next we consider the case of a non-constant $\epsilon_1|_{\R_+}$ or $\epsilon_1|_{\R_-}$. Numerically we find
examples for which \eqref{eigenvalue_problem_numerics} possesses a localised solution, again provided that a suitable (non-explicit) dispersion relation $k=k(\omega)$ is fulfilled.

\paragraph{Numerical implementation} To simplify the numerics we write $w_3$ as the sum of a $C^1$ component $w_{3,r}$ and a component $w_{3,s}$ of a simple form and with a discontinuous derivative. $w_{3,s}$ can be expressed in term of $w_{3,r}$ and hence we only need to solve for $w_{3,r}$. In detail, let $w_3=w_{3,r}+w_{3,s}$ with
\begin{equation*}
w_{3,s}(x_1) = \begin{cases}
w_{3,s}^-\in\R & x_1 < 0,\\
w_{3,s}^+(x_1) & x_1 \geq0,
\end{cases}
\end{equation*}
and choose the constant $w_{3,s}^-$ such that $w_{3,s}$ is continuous, i.e. $w_{3,s}^- =w_{3,s}^+(0)$. Note that with this choice $w_3$ is continuous. Moreover,
\begin{align}
\left\llbracket \frac{w_3'}{\epsilon_1} \right\rrbracket = 0\Leftrightarrow &\quad\epsilon_1^-(0)\left(w_{3,r}'(0) + (w_{3,s}^+)'(0)\right) = \epsilon_1^+(0)\,w_{3,r}'(0) \notag\\
\Leftrightarrow &\quad(w_{3,s}^+)'(0) = \frac{\epsilon_1^+(0) - \epsilon_1^-(0)}{\epsilon_1^-(0)}\,w_{3,r}'(0) =:\nu\,w_{3,r}'(0).\label{E:w3sr}
\end{align}
If we set
\begin{equation*}
w_{3,s}(x_1)=\left(\mathcal{L}w_{3,r}\right)(x_1):= \begin{cases}
- \sign{\nu}\,w_{3,r}'(0) & x_1 < 0,\\
- \sign{\nu}\,w_{3,r}'(0)\,\e^{-|\nu| x_1} & x_1 \geq0,
\end{cases}
\end{equation*}  
then \eqref{E:w3sr} holds and the second interface condition in \eqref{eigenvalue_problem_numerics} is satisfied. Note that $\cL:C^1(\R)\to C(\R)$ is a \textit{linear} operator and the equation for $w_{3,r}$ is the differential equation
\begin{equation}\label{eigenvalue_problem_numerics_2}
\begin{cases}
\left(-\partial_{x_1}^2 + \frac{\epsilon_1'}{\epsilon_1} \partial_{x_1} - \epsilon_1\omega^2 \right)\left((I + \mathcal{L})w_{3,r}\right)=-k^2 (I + \mathcal{L})w_{3,r}&\mbox{in}\quad\R\setminus\{0\},\\
\llbracket w_{3,r}\rrbracket=\left\llbracket w_{3,r}'\right\rrbracket=0.%
\end{cases}
\end{equation}
We look for solutions $w_3$ which are eigenfunctions, i.e. with $w_3\in H^1(\R)$. In the special case when $w_3(x_1)\to0$ for $|x_1|\to+\infty$, the corresponding conditions for $w_{3,r}$ are
\begin{equation*}
\lim_{x_1\to-\infty}w_{3,r}(x_1)=\sign{\nu}\,w_{3,r}'(0),\qquad\lim_{x_1\to+\infty}w_{3,r}(x_1)=0.
\end{equation*}
Note that we can freely choose $w_{3,r}'(0)$.

To solve \eqref{eigenvalue_problem_numerics_2} numerically for a fixed $\omega\in\R$ we discretise the problem and apply a solver for a generalised eigenvalue problem, e.g. a solver based on a Krylov-Schur algorithm.

\begin{figure}[ht!]
\captionsetup{font=footnotesize}
	\centering
	\begin{subfigure}[b]{0.38\linewidth}
		\includegraphics[width=\linewidth]{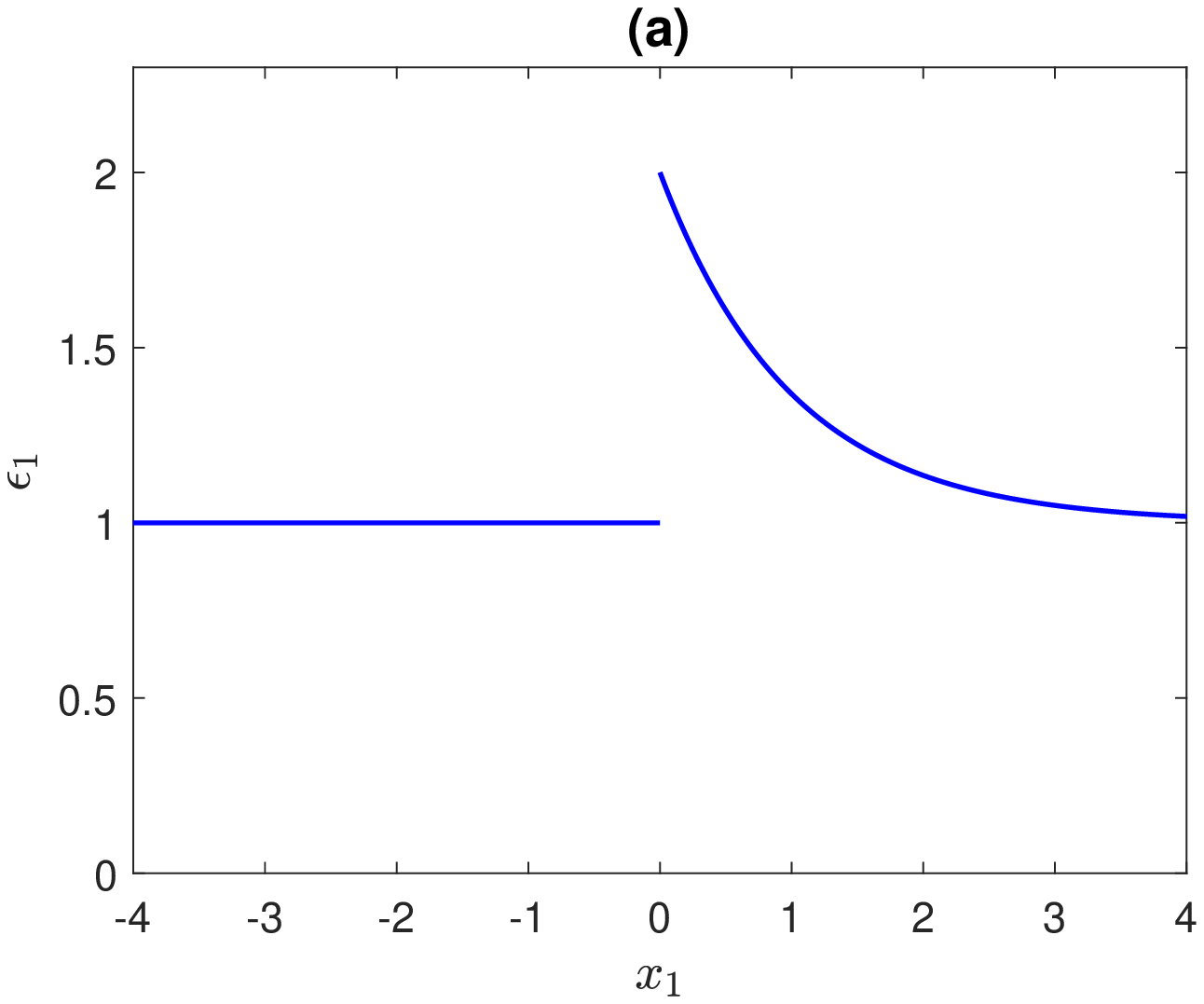}
	\end{subfigure}
	$\;$
	
	\begin{subfigure}[b]{0.38\linewidth}
		\includegraphics[width=\linewidth]{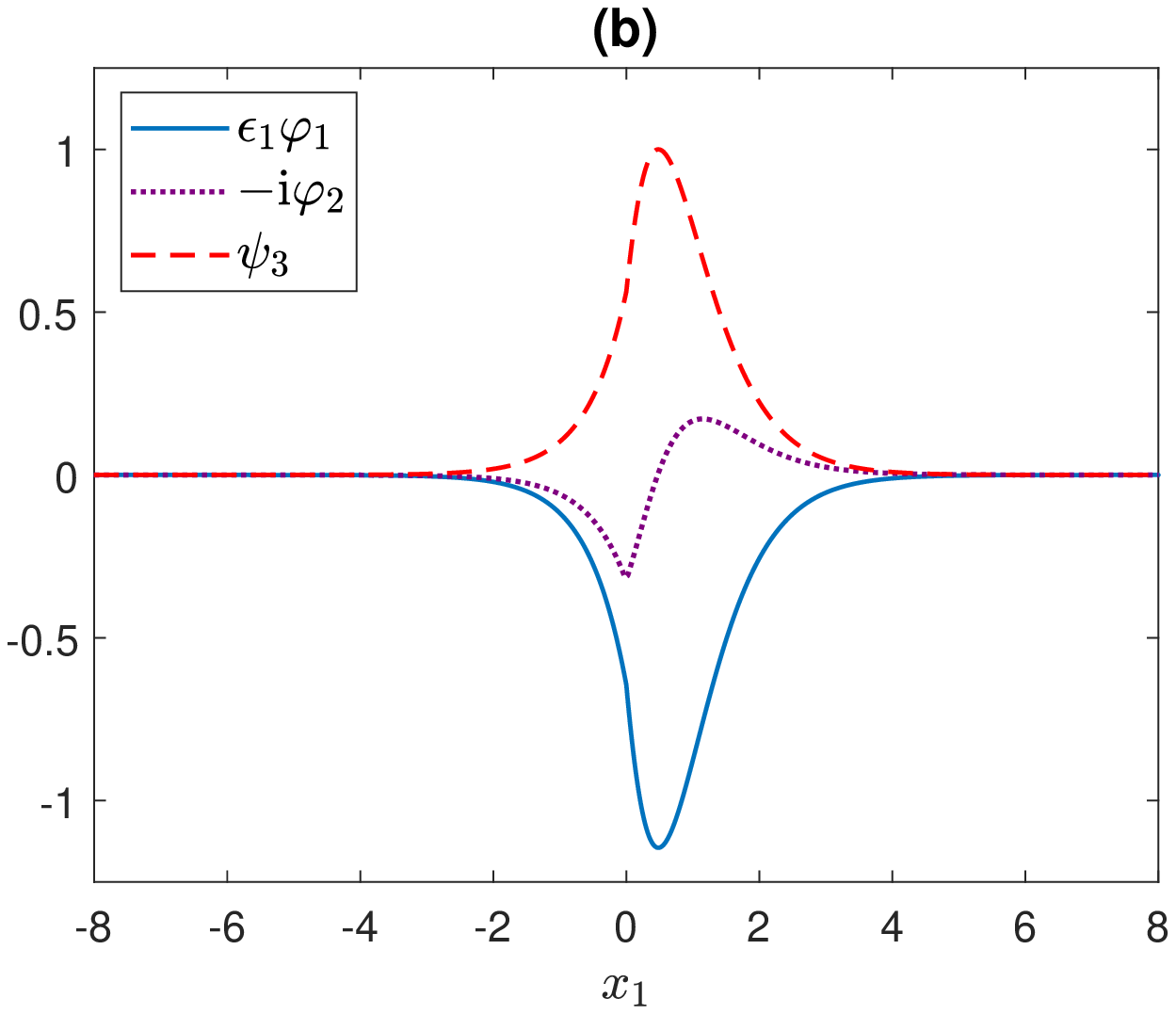}
	\end{subfigure}
	$\;$
	\begin{subfigure}[b]{0.38\linewidth}
		\includegraphics[width=\linewidth]{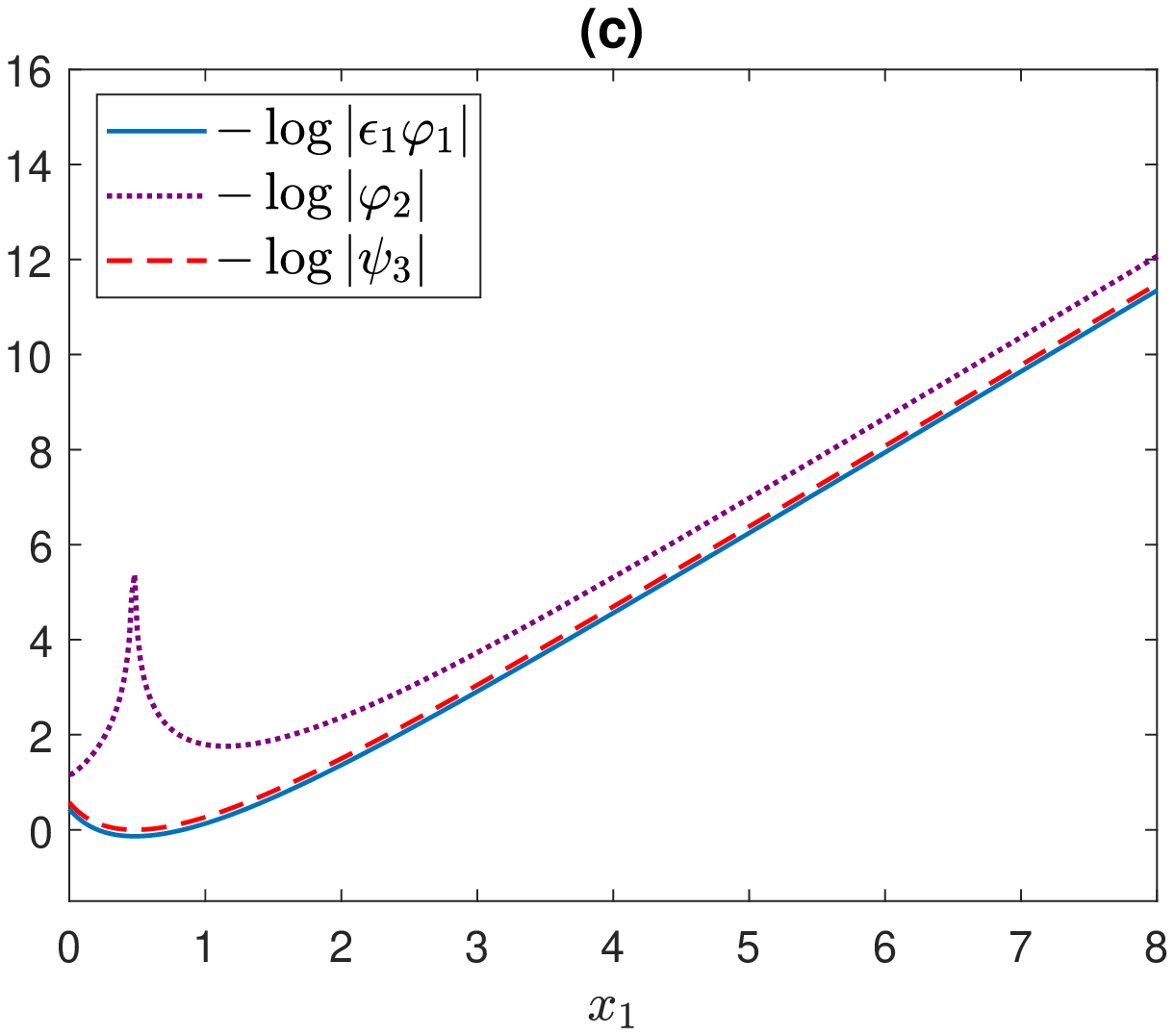}
	\end{subfigure}
	\caption{ (a) The dielectric function $\epsilon_1(x_1)=1\chi_{\R_-}+\left(1+{\rm e}^{-x_1}\right)\chi_{\R_+}$. (b) the respective solution $w(x_1)=\left(\varphi_1(x_1),\varphi_2(x_1),\psi_3(x_1)\right)^\trans$ of the linear problem \eqref{L_selfadj_eq} for $\omega_0=3$ (cf. Remark \ref{real_mag_real}). A corresponding eigenvalue $k_0=k_0(\omega_0)$ produced by the numerics is $k_0\approx3.4352$. We plot $\epsilon_1\varphi_1$ instead of $\varphi_1$ to show the fact all interface conditions in \eqref{L_selfadj_eq} are satisfied. (c) Demonstration of the exponential decay of $\epsilon_1\varphi_1$, $\varphi_2$, and $\varphi_3$ for $x_1\to+\infty$.}
	\label{Fig_localised}
\end{figure}

In Figure \ref{Fig_localised}(a) we plot a possible choice of a dielectric function $\epsilon_1$ satisfying assumptions (H0)-(H1) and with a jump at the interface. We chose $\epsilon_1$ constant for $x_1<0$ and of the form $1+\e^{-x_1}$ for $x_1\geq0$. The respective solution of \eqref{L_selfadj_eq} with an arbitrary normalisation is plotted in Figure \ref{Fig_localised}(b) and is exponentially decaying at $\pm\infty$, as shown in (c) for $x_1>0$. Note that it is a simple exercise to prove that the solution is exponentially decaying for $x_1\to-\infty$ since here the potential is constant and the calculations may be carried out explicitly. We also point out that the derivative of $\varphi_2$ has a jump at $x_1=0$. This is due to the fact that by the third equation \eqref{L_selfadj_eq} one has $-\ri\varphi_2'=k\varphi_1+\omega\psi_3$, a sum of a discontinuous and a continuous quantity.

\subsection{Application of Theorems \ref{Thm_phi}-\ref{Thm_der_phi}}\label{Section_ApplicationNL}
Let us choose $\omega_0\in\R$ and $k_0=k_0(\omega)$ as one of the corresponding eigenvalues. The corresponding eigenfunction as introduced in Section \ref{Section_linear} is $\left(\varphi_1,\varphi_2,\psi_3\right)^\trans$. By our choice of TM-modes, the vector $m_\cU$ in the ansatz \eqref{Uans0} is defined as $m_\cU(x_1):=\left(\varphi_1(x_1),\varphi_2(x_1),0,0,0,\psi_3(x_1)\right)^\trans$ and, dropping the third vanishing component, we may redefine the vector $m_\cE$ in \eqref{Eans0} as $m_\cE(x_1):=\left(\varphi_1(x_1),\varphi_2(x_1)\right)^\trans$. Hence our given vector field $U_0=\left(U_{0,1},U_{0,2}\right)^\trans:\R^2\to\R^2$, with which one begins the nonlinear analysis, is
\begin{equation*}
U_0(x_1,x_2)=\varepsilon\cA(\varepsilon x_2)m_\cE(x_1)\e^{\ri k_0x_2}+\text{c.c.}\,,
\end{equation*}
where $U_0$ needs to satisfy assumptions of Theorems \ref{Thm_phi}-\ref{Thm_der_phi}, which is the goal of this section. In the amplitude approximation of the wavepacket \eqref{Uans0}, $\cA$ is the initial condition $A(\cdot,0)$. We do not study the effective Schr\"odinger equation for $A$ in this paper, see \cite{KSM,DR,MS,LS}.

We underline the fact that $U_0$ scales differently in $x_1$ and $x_2$ with respect to $\varepsilon$. It is easy to see that
\begin{equation*}
\|U_0\|_2=\cO(\varepsilon^{1/2})
\end{equation*}
provided $m_\cE\in L^2(\R)^2$ and $\cA\in L^2(\R)$. The loss of the $\tfrac12$-power of $\varepsilon$ is due to the scaling in the second variable. Hence, in order for the solution $\nabla\phi$ of \eqref{eq}-\eqref{IFCs_eq} to be a correction term of $U_0$, its $L^2$-norm should be $o(\varepsilon^{1/2})$. The reason why we may expect such a behaviour for $\nabla\phi$ is that, if we isolate the term $b=\nabla\cdot\left(\epsilon_1 U_0\right)$ as in \eqref{eq}, we get
\begin{equation*}
\begin{split}
b(x_1,x_2)&= \partial_1 \left(\epsilon_1\Aou\right) + \epsilon_1 \partial_2\Aod\\
&=\left[\varepsilon\cA(\varepsilon x_2)\left( \left(\epsilon_1 m_{\cE,1}\right)' + \epsilon_1 \ri k_0 m_{\cE,2}\right)+\varepsilon^2\epsilon_1\cA'(\varepsilon x_2)m_{\cE,2}\right]\e^{\ri k_0x_2}+\text{c.c.}\,.
\end{split}
\end{equation*}
In the first term we may recognise $\nabla\cdot\cD_{lin}$ as in \eqref{divDlin=0}, so this term vanishes and hence
\begin{equation}\label{b_epsilon_2}
b(x_1,x_2)=\varepsilon^2\epsilon_1(x_1)\cA'(\varepsilon x_2)m_{\cE,2}(x_1)\e^{\ri k_0x_2}+\text{c.c.}\,,
\end{equation}
which implies
\begin{equation}\label{b_2}
\|b\|_2=\cO(\varepsilon^{3/2}),
\end{equation}
provided $m_{\cE,2}\in L^2(\R)$ and $\cA'\in L^2(\R)$. Notice that if we were in a smooth bounded domain, e.g. with Dirichlet boundary conditions, then the estimate $\|\phi\|_{H^2}\leq C\varepsilon^{3/2}$ would follow from the weak formulation \eqref{wf} directly using Poincar\'{e} inequality. On the other hand, if one looks at the weak formulation \eqref{wf} with $\eta=\phi$ and naively tries to estimate $\|\nabla\phi\|_2^2$ by the absolute value of the right-hand side and then using the Cauchy-Schwarz inequality, one ends up only with the unsatisfactory estimate $\|\nabla\phi\|_2\leq C\varepsilon^{1/2}$, because $\|\epsilon_1U_0\|_2=\cO(\varepsilon^{1/2})$. This is the reason why it is important that the right-hand sides of the estimates of Theorems \ref{Thm_phi}-\ref{Thm_der_phi} are not dependent on $\|U_0\|_2$, as briefly observed in Remark \ref{Remark_smallness}. With our choice of a cubic $\cD$-field as in \eqref{Kerr}, one may first apply Theorem \ref{Thm_phi} with $\alpha=p=3$ (since there are no further nonlinear terms in $\cD$) and get, in particular,
\begin{equation}\label{nablaphi_2_Anwendung}
\|\nabla\phi\|_2^2\les\|U_0\|_4^4+\|b\|_2^2+\|b\|_{L^1(\log)}^2,
\end{equation}
where we recall that $\|b\|_{L^1(\log)}=\int|b(x)|\log(2+|x|)\,dx$. The term $\|b\|_2$ was estimated in \eqref{b_2}. Next,
\begin{equation}
\int|U_0|^4=\varepsilon^4\int_\R|\cA(\varepsilon x_2)|^4\,dx_2\int_\R|m_\cE(x_1)|^4\,dx_1\leq\varepsilon^3\|\cA\|_4^4\|m_\cE\|_4^4,
\end{equation}
so provided $\cA\in L^4(\R)$ and $m_\cE\in L^4(\R)^2$, one infers
\begin{equation}\label{A0_4}
\|U_0\|_4=\cO(\varepsilon^{3/4}).
\end{equation}
For the last term in \eqref{nablaphi_2_Anwendung}, we first note the simple equality
$$\log(s+t)=\log(s)+\log\Big(1+\frac ts\Big),\quad\;\,s,t\in\R^+,$$
to obtain, for a suitable constant $c>0$,
\begin{equation*}
\begin{split}
\log(2+|x|)&\leq\log(2+c|x_1|+c|x_2|)=\log(2+c|x_1|)+\log\left(1+\frac{c|x_2|}{2+c|x_1|}\right)\\
&\leq\log(2+c|x_1|)+\log(2+c|x_2|).
\end{split}
\end{equation*}
Hence,
\begin{equation*}
\begin{split}
\|b\|_{L^1(\log)}&\leq\varepsilon^2\|\epsilon_1\|_\infty\int\left|\cA'(\varepsilon x_2)\right||m_{\cE,2}(x_1)|\log(2+|x|)\,dx\\
&\les\varepsilon^2\|m_{\cE,2}\|_{L^1(\log)}\int\left|\cA'(\varepsilon x_2)\right|\,dx_2+\varepsilon^2\|m_{\cE,2}\|_1\int\left|\cA'(\varepsilon x_2)\right|\log(2+c|x_2|)\,dx_2\\
&\les\varepsilon^2\|m_{\cE,2}\|_{L^1(\log)}\int\left|\cA'(\varepsilon x_2)\right|\left(1+\log(2+c|x_2|)\right)\,dx_2\\
&\les\varepsilon\|m_{\cE,2}\|_{L^1(\log)}\int\left|\cA'(y)\right|\log\left(e+c\frac{|y|}\varepsilon\right)\,dy,
\end{split}
\end{equation*}
where in the last inequality we made use of the variable transformation $y=\varepsilon x_2$. Finally, choosing any $\gamma>0$ one may estimate $\log\left(e+ct\right)\les1+t^\gamma$ for all $t\geq0$. Therefore
\begin{equation}\label{b_log}
\begin{split}
\|b\|_{L^1(\log)}&\les\|m_{\cE,2}\|_{L^1(\log)}\left(\varepsilon\|\cA'\|_1+\varepsilon^{1-\gamma}\int\left|\cA'(y)\right||y|^\gamma\,dy\right)=\cO(\varepsilon^{1-\gamma}),
\end{split}
\end{equation}
provided $m_{\cE,2}\in L^1(\log,\R)$ and $\cA'\in L^1(\gamma,\R):=\{\varphi\in L^1(\R)\,|\,\int|\varphi(x)|\left(1+|x|^\gamma\right)\,dx<\infty\}$. In summary, from estimate \eqref{nablaphi_2_Anwendung}, one infers by \eqref{b_2}, \eqref{A0_4}, and \eqref{b_log} we infer
\begin{equation*}
\|\nabla\phi\|_2\leq C_1\varepsilon^{3/2}+C_2\varepsilon^{1-\gamma},
\end{equation*}
where $\gamma>0$ is arbitrary and the constants $C_1$ and $C_2$, do not depend on $\varepsilon$. Hence, choosing $\gamma\in\left(0,\frac12\right)$ one obtains the desired estimate of order $o(\varepsilon^{1/2})$. For $\gamma\in (0,\tfrac12)$ we get $\|\nabla\phi\|_2=O(\varepsilon^{1-\gamma})$.

Analogous estimates may be deduced for the first derivatives of $\nabla\phi$ according to Theorem \ref{Thm_der_phi} if one prescribes further regularity on the functions $\cA$ and $m_\cE$. For the sake of clarity, we collect the assumptions needed to estimate the quantities $\|U_0\|_4$, $\|b\|_2$, $\|b\|_{L^1(\log)}$ in Table \ref{table_regularity_assumptions}. For $n,m\in\N\setminus\{0\}$ we denote $W^{1,p}_\pm(\R^n)^m:=\{\varphi\in L^p(\R^n)^m\,|\,\partial_k \varphi \in L^p(\R^n)^m\;\mbox{for}\;k\geq2,\,\partial_1\varphi\in L^p(\R^n_\pm)^m\}$ and define the norm accordingly as $\|\varphi\|_{W^{1,p}_\pm}:=\sum_{k=2}^n\|\partial_k\varphi\|_p+\sum_{\pm}\|\partial_1\varphi\|_{L^p(\R^n_\pm)}$.

\begin{table}[H]
	\captionsetup{font=footnotesize}
	\centering
	\begin{tabular}{ccccc}
		\toprule
		$\|U_0\|_4$ & $\|b\|_2$ & $\|b\|_{\log}$ & $\|U_0\|_{W^{1,4}_\pm}$ & $\|\partial_2b\|_2$\\
		\midrule
		$\cA\in L^4$ & $\cA'\in L^2$ & $\cA'\in L^1(\gamma)$ & $\cA\in W^{1,4}$ &$\cA'\in H^1$\\
		$m_\cE\in L^4$ & $m_{\cE,2}\in L^2$ & $m_{\cE,2}\in L^1(\log)$ & $m_\cE\in W^{1,4}_\pm$ & $m_{\cE,2}\in L^2$\\
		\bottomrule
	\end{tabular}
	\caption{Sufficient regularity assumptions on $\cA$ and $m_{\cE}$ to estimate the terms on the first line.}
	\label{table_regularity_assumptions}
\end{table}

Summarising our results, we are led to the following.

\begin{prop}\label{Prop_Anwendung}
Assume (H0)-(A4) and let $\cA$ and $m_\cE$ satisfy all the conditions in Table \ref{table_regularity_assumptions}. Then there exists a function $\phi\in\fD_{2,4}$ such that $V = \nabla \phi$ is a solution of the problem \eqref{E0}-\eqref{initial_div}-\eqref{Max_IFCs_korr} for $n=2$ with $\|V\|_{W^{1,2}_\pm}\leq C\varepsilon^{1-\gamma}$ for any $\gamma>0$, where \eqref{initial_div} and the third interface condition in \eqref{Max_IFCs_korr} hold in the sense of Definitions \ref{weak_sol} and  \ref{T+=T-} resp. and the remaining interface conditions hold pointwise.
\end{prop}
\begin{remark}Due to $\|\nabla\phi\|_{L^2}=o(\varepsilon^{1/2})$ the function $\nabla \phi$ may be interpreted as a correction of the initial value $U_0=\cE_{ans}^{(0)}$.
\end{remark}
\begin{proof}
	It just remains to show that $\llbracket\cE^{(0)}_2\rrbracket=0$ since the third component of the electric field vanishes for TM-modes. This condition is already satisfied by $U_0$ because of the choice of $m_\cE$, so one needs show that $V_2 = \partial_2\phi$ is continuous across the interface. However, since $\phi\in\fD_{2,4}\subset D^{1,4}_0(\R^2)$, by \eqref{estimateLpD1p} we know that $\phi\in W^{1,4}_{loc}(\R^2)$, so in particular $\phi\in C^{0,\frac12}_{loc}(\R^2)$ by the Morrey embedding. This means that $\phi$ is continuous across the interface and, further, that every tangential derivative with respect to the interface is continuous too, see e.g. \cite[Sec 173-175]{TT}.
\end{proof}

\begin{remark}
It clearly follows from the analysis above and Remark \ref{Remark_higher_regularity} that it is possible to obtain analogous bounds with the same order in $\varepsilon$ for the $H^k$-norm of $\nabla\phi$, provided higher regularity of $\cA$ in $\R$ and $m_\cE$ in $\R_\pm$ is prescribed.
\end{remark}

\begin{remark}
	Note that the exponential decay of $m_\cE$ shown by the numerics (Figure \ref{Fig_localised}) indicates that the assumptions on $m_\cE$ in Table \ref{table_regularity_assumptions} are satisfied. In the application where $\cA$ is the initial condition of the envelope $A$, and hence free to choose, the conditions on $\cA$ in Table \ref{table_regularity_assumptions} are not critical.
\end{remark}

\begin{remark} \label{Remark_3/2}
	Continuing the parallel with the case of a bounded domain, in which it is clear that $\|b\|_2=\cO(\varepsilon^{3/2})$ would imply $\|\nabla\phi\|_{H^1}=\cO(\varepsilon^{3/2})$, and hence $\nabla\phi$ is of the same order in $\varepsilon$ as the right-hand side, we expect that the estimates provided by Theorems \ref{Thm_phi}-\ref{Thm_der_phi} are not optimal. In our case $\alpha=p=3$ we can see that this is due to the additional term $\iii b$, which produces the logarithmic term as $n=2$, while all other terms on the right-hand side are indeed of order $\cO(\varepsilon^{3/2})$. Thus it would be desirable to be able to estimate the term $b$ differently than by Proposition \ref{OS_thm}.
\end{remark}

\subsection{Numerical test}\label{Section_numericaltest}
To confirm the analytic results of Propositions \ref{Prop_Anwendung}, we calculate a numerical solution of problem \eqref{E0}-\eqref{initial_div}-\eqref{Max_IFCs_korr} by applying a fixed point iteration to a finite element discretisation of the problem.
 
For a cubic nonlinearity, as described in the beginning of this chapter, we therefore have to solve
\begin{equation}
\begin{cases}
- \nabla \cdot \left(\epsilon_1 \nabla \phi \right) = \mathfrak{f}(\phi) & \mbox{in}\quad\R_\pm^2, \\
\left\llbracket \epsilon_1 \partial_{1} \phi \right\rrbracket = \mathfrak{h}(\phi), \\
\left\llbracket \partial_{2} \phi \right\rrbracket = 0,
\end{cases}
\label{interface_problem_numeric}
\end{equation}
with the $\phi$-dependent functions
\begin{align*}
\mathfrak{f}(\phi) &:= \nabla \cdot \left(\epsilon_1 U_0 + \epsilon_3 \left|U_0 + \nabla \phi\right|^2 \left(U_0 + \nabla \phi \right) \right),\\
\mathfrak{h}(\phi) &:= - \left\llbracket \epsilon_3 \left|U_0 + \nabla \phi\right|^2 \left(U_{0,1} + \partial_{1} \phi \right)\right\rrbracket.
\end{align*}
To find a solution of \eqref{interface_problem_numeric}, we rewrite the problem as a system of two coupled Neumann boundary value problems, in which we have to determine the functions $\phi: \R^2_\pm \rightarrow \R$ and $\mathfrak{g}: \Gamma \rightarrow \R$ such that
\begin{equation}
\begin{cases}
- \nabla \cdot \left(\epsilon_1 \nabla \phi \right) = \mathfrak{f}(\phi) & \mbox{in}\quad\R_\pm^2, \\
\left(\epsilon_1 \partial_1 \phi\right)_- = \mathfrak{g}, \\
\left(\epsilon_1 \partial_1 \phi\right)_+ = \mathfrak{h}(\phi) + \mathfrak{g}, \\
\left\llbracket \partial_{2} \phi \right\rrbracket = 0,
\end{cases}
\label{interface_problem_numeric_2}
\end{equation}
where we recall $u_\pm(x):=\lim_{t\to 0 \pm}u(x + te_1)$, for $x \in \Gamma$. Note that a solution $\phi$ of \eqref{interface_problem_numeric_2} is also a solution of \eqref{interface_problem_numeric}.

We now approximate the solution of the nonlinear problem \eqref{interface_problem_numeric_2} with the help of a fixed point iteration. We select an initial guess $\phi_0$ and solve
\begin{equation}
\begin{cases}
- \nabla \cdot \left(\epsilon_1 \nabla \phi_{n+1} \right) = \mathfrak{f}(\phi_n) & \mbox{in}\quad\R_\pm^2, \\
\left(\epsilon_1 \partial_1 \phi_{n+1}\right)_- = \mathfrak{g}_{n+1}, \\
\left(\epsilon_1 \partial_1 \phi_{n+1}\right)_+ = \mathfrak{h}(\phi_n) + \mathfrak{g}_{n+1}, \\
\left\llbracket \partial_{2} \phi_{n+1} \right\rrbracket = 0
\end{cases}
\label{interface_problem_numeric_fixed_point}
\end{equation}
iteratively for $n \geq 0$. The weak formulation of the problem is given by 
\begin{equation}
\begin{cases}
\int_{\R^2_+} \epsilon_1 \nabla \phi_{n+1} \cdot \nabla \eta + \int_\Gamma \mathfrak{g}_{n+1} \eta = \int_{\R^2_+} \mathfrak{f}(\phi_n) \eta - \int_\Gamma \mathfrak{h}(\phi_n) \eta & \eta \in H^1(\R^2_+), \\
\int_{\R^2_-} \epsilon_1 \nabla \phi_{n+1} \cdot \nabla \eta - \int_\Gamma \mathfrak{g}_{n+1} \eta = \int_{\R^2_-} \mathfrak{f}(\phi_n) \eta & \eta \in H^1(\R^2_-), \\
\left\llbracket \partial_{2} \phi_{n+1} \right\rrbracket = 0.
\end{cases}
\label{interface_problem_numeric_weak}
\end{equation}
To solve \eqref{interface_problem_numeric_weak} numerically, we use the finite element method. First we replace $\R^2_\pm$ and $\Gamma$ by suitable bounded domains $\Omega_\pm \subset {\overline{\R^2_\pm}}$ and $\widetilde{\Gamma} := \overline{\Omega}_+ \cap \overline{\Omega}_- \subset \Gamma$, respectively. Furthermore we substitute $H^1(\R^2_\pm)$ with the following $N$-dimensional subspaces $\cV_\pm := \operatorname{span}\{\eta_k^\pm\,|\, k = 1,\dots,N\}$, $N \in \N$, where the shape functions $\eta_k^\pm \in H^1(\Omega_\pm)$ are the standard piecewise linear hat functions, which are linearly independent. Then we look for solutions of the form
\begin{align*}
\phi_{n+1}(x) &= \begin{cases}
\sum_k \Phi_{n+1,k}^+\,\eta_k^+(x), & x \in \Omega_+,\\
\sum_k \Phi_{n+1,k}^-\,\eta_k^-(x), & x \in \Omega_-,\\
\end{cases} \\
\mathfrak{g}_{n+1}(x) &= \sum_k G_{n+1,k}\,\eta_k^+(x)|_{\widetilde{\Gamma}},
\end{align*} 
where the coefficients $\Phi_{n+1,k}^\pm, G_{n+1,k}$ are the solutions of the following system
\begin{equation*}
\begin{cases}
\sum_k \left( \Phi_{n+1,k}^+ \int_{\Omega_+} \epsilon_1 \nabla \eta_j^+ \cdot \nabla \eta_k^+ + G_{n+1,k} \int_{\widetilde{\Gamma}} \eta_j^+ \eta_k^+ \right) = \int_{\Omega_+} \mathfrak{f}(\phi_n) \eta_j^+ - \int_{\widetilde{\Gamma}} \mathfrak{h}(\phi_n) \eta_j^+ & \quad j=1,\dots,N, \\
\sum_k \left( \Phi_{n+1,k}^- \int_{\Omega_-} \epsilon_1 \nabla \eta_j^- \cdot \nabla \eta_k^- - G_{n+1,k} \int_{\widetilde{\Gamma}} \eta_j^- \eta_k^+ \right) = \int_{\Omega_-} \mathfrak{f}(\phi_n) \eta_j^- & \quad j=1,\dots,N, \\
\sum_k \left( \Phi_{n+1,k}^- \partial_2 \eta_k^-|_{\widetilde{\Gamma}} - \Phi_{n+1,k}^+ \partial_2 \eta_k^+|_{\widetilde{\Gamma}} \right)= 0.
\end{cases}
\end{equation*}  
Due to the Neumann boundary conditions and the fact that $\mathfrak{f}$ and $\mathfrak{h}$ depend only on the gradient of $\phi$, the solution is unique only up to an additive constant. To get uniqueness we additionally demand $\int_{\Omega_\pm} \phi_{n+1} = 0$. Therefore, we extend the finite element formulation by 
\begin{align*}
\sum_k \Phi_{n+1,k}^\pm \int_{\Omega_\pm} \eta_k^\pm = 0.
\end{align*}  

Let us now consider the example introduced in Section \ref{Section_linear} with the main goal of verifying Proposition \ref{Prop_Anwendung} and in particular the conjectured convergence rate $\mathcal{O}(\varepsilon^{3/2})$ for $\|\nabla \phi\|_2$  in Remark \ref{Remark_3/2}. Due to the fact that we choose piecewise linear shape functions we can only study the convergence of the $L^2$-norm of $\nabla\phi$ although Proposition \ref{Prop_Anwendung} estimates even the $H^1$-norm.
For our numerical test we choose $\epsilon_1(x_1)=1\chi_{\R_-}+\left(1+{\rm e}^{-x_1}\right)\chi_{\R_+}$ and $m_\mathcal{E}$ accordingly, see Figure \ref{Fig_localised}, $\epsilon_3  \equiv 1$, and we set $\mathcal{A}(x_2) = \e^{-5\cdot 10^6 \varepsilon^2 x_2^2}$. Note that with such choices of (exponentially decaying) $m_{\mathcal{E}}$ and $\mathcal{A}$, the conditions for $U_0$ and  $b$ in Table \ref{table_regularity_assumptions} are satisfied. For the discretisation we select $\overline{\Omega}_- = [-6,0] \times [-6,6]$ and $\overline{\Omega}_+ = [0,6] \times [-6,6]$ and choose a regular triangulation of step size $h$ together with standard hat-functions for $\eta_k^\pm$. At those boundaries of $\Omega_\pm$ which are not part of $\widetilde{\Gamma}$ we enforce homogeneous Neumann boundary conditions. 

\begin{figure}[ht!]
	\captionsetup{font=footnotesize}
	\centering
	\begin{subfigure}[b]{0.38\linewidth}
		\includegraphics[width=\linewidth]{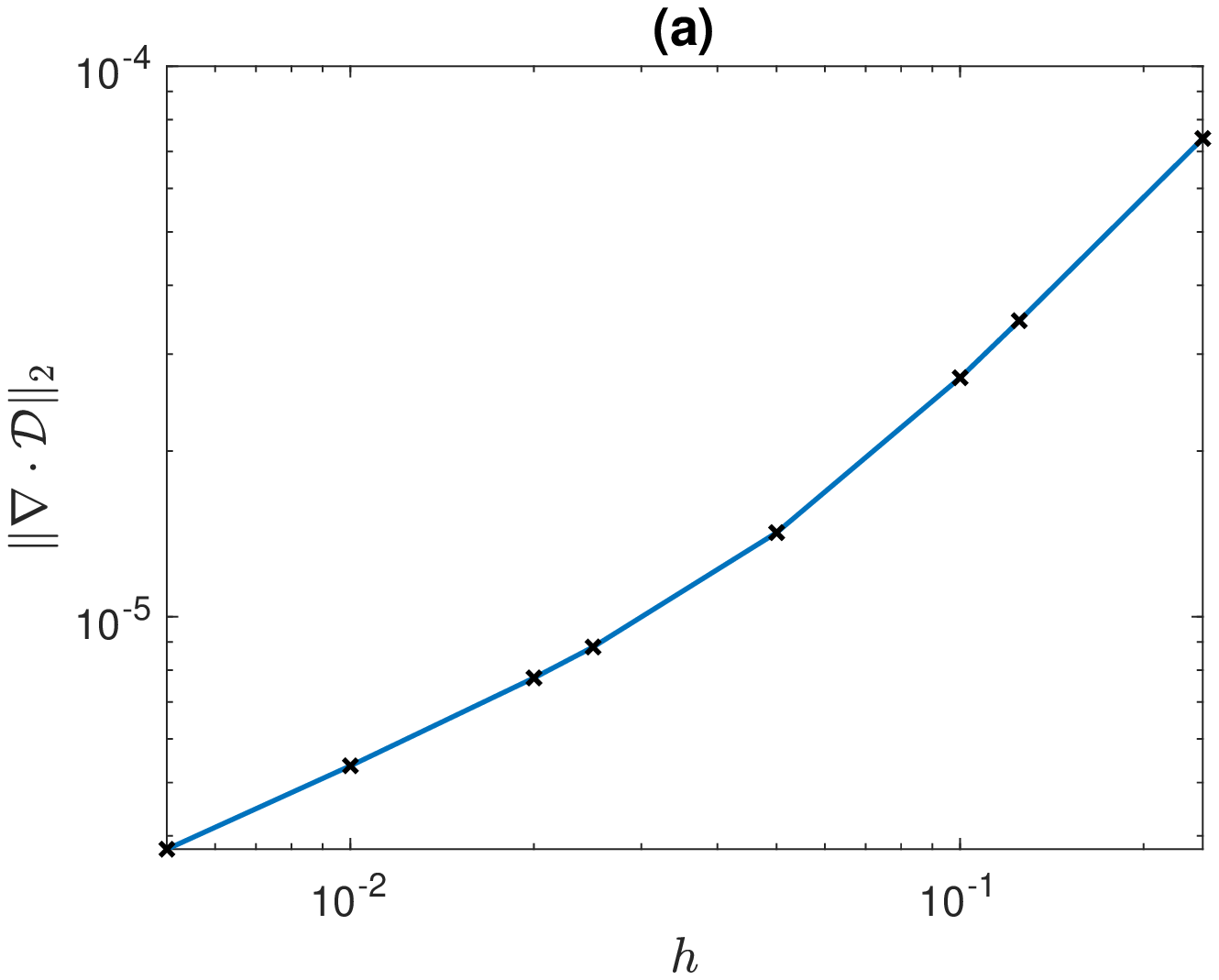}
	\end{subfigure}
	$\;$
	\begin{subfigure}[b]{0.38\linewidth}
		\includegraphics[width=\linewidth]{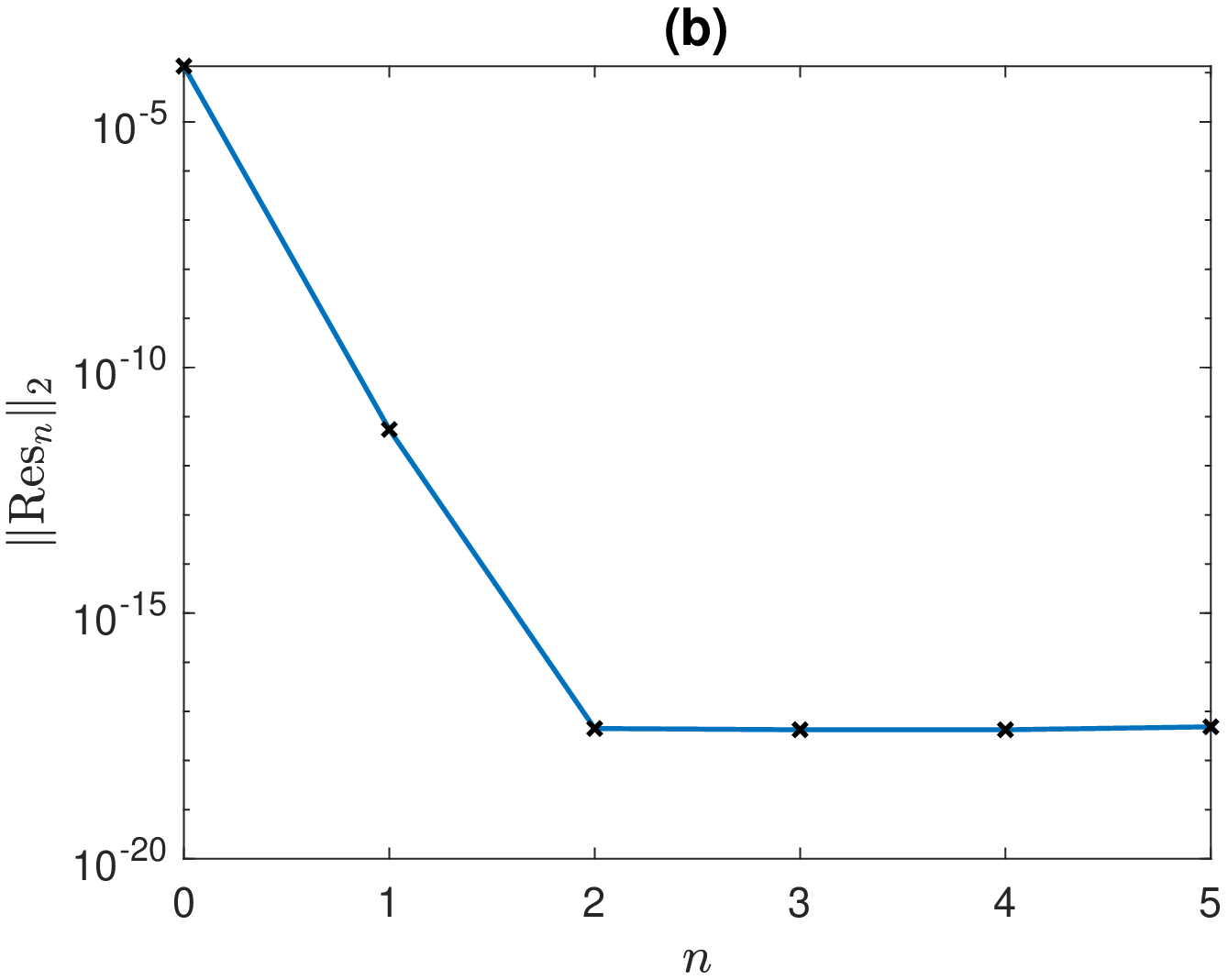}
	\end{subfigure}
	\caption{(a) $\| \nabla \cdot \mathcal{D}\|_2$ in dependence on the step size $h$.  (b) Plot of the $L^2$-norm of the residual in the fixed point iteration for each step.}
	\label{Fig_fixed_point}
\end{figure}

For the fixed point iteration we start with $\phi_0 \equiv 0$ as an initial guess. Let us first check the convergence of the discretisation (in $h$ and in the iteration $n$). For Figure \ref{Fig_fixed_point} (a) we fixed $\varepsilon = 3 \cdot 10^{-4}$ and calculated $\| \nabla \cdot \mathcal{D}\|_2$ for different step sizes $h$ ranging from $0{.}25$ to $0{.}005$. For Figure \ref{Fig_fixed_point} (b) we also fixed $h=0{.}005$ and calculated the $L^2$-norm of the residual $\mathrm{Res}_n := - \nabla \cdot \left(\epsilon_1 \nabla \phi_{n} \right) - \mathfrak{f}(\phi_n)$ in each step of the fixed point iteration. We see the numerical convergence in both plots. \\

\begin{figure}[ht!]
\captionsetup{font=footnotesize}
	\centering
	\begin{subfigure}[b]{0.38\linewidth}
		\includegraphics[width=\linewidth]{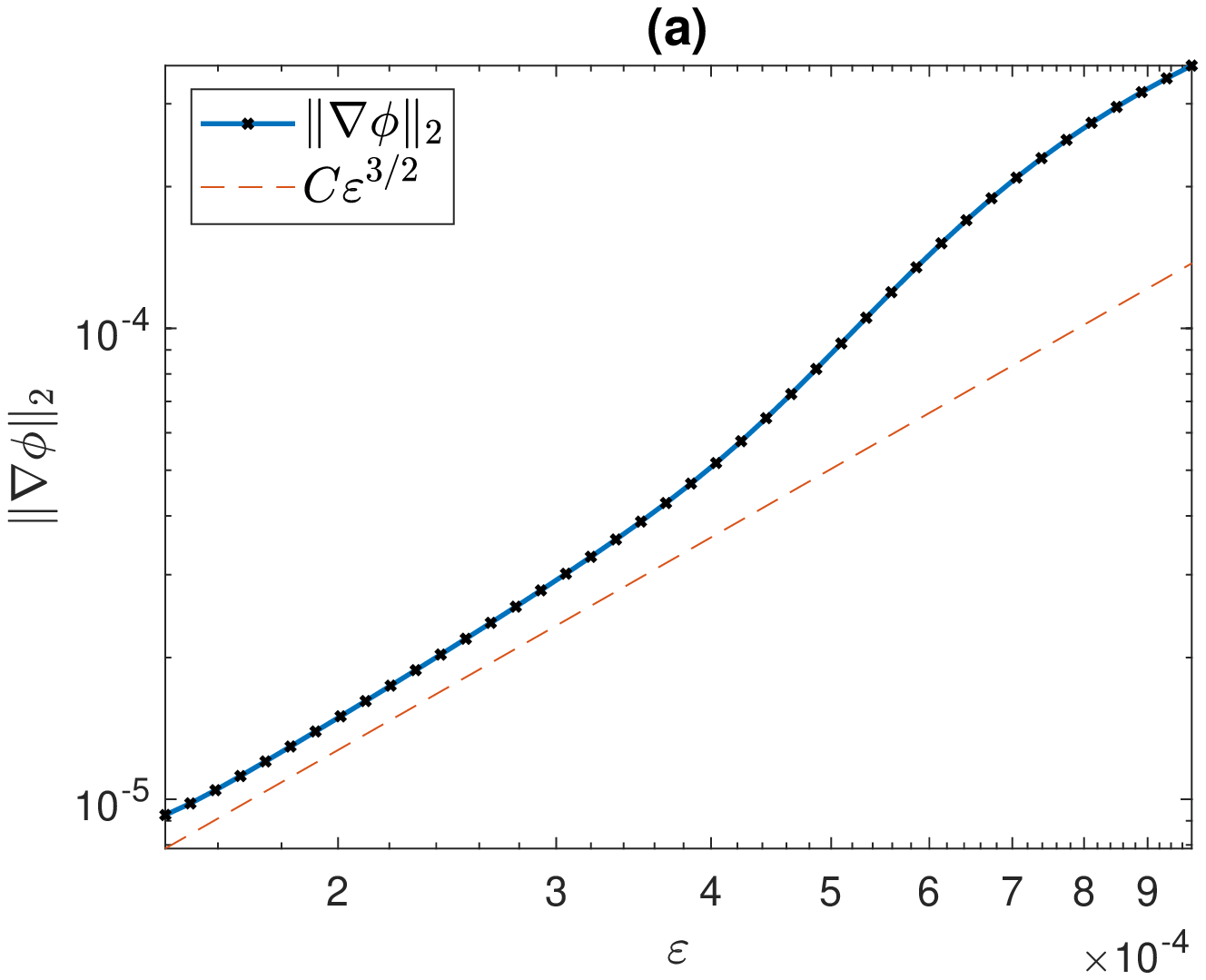}
	\end{subfigure}
	$\;$
	
	\begin{subfigure}[b]{0.38\linewidth}
		\includegraphics[width=\linewidth]{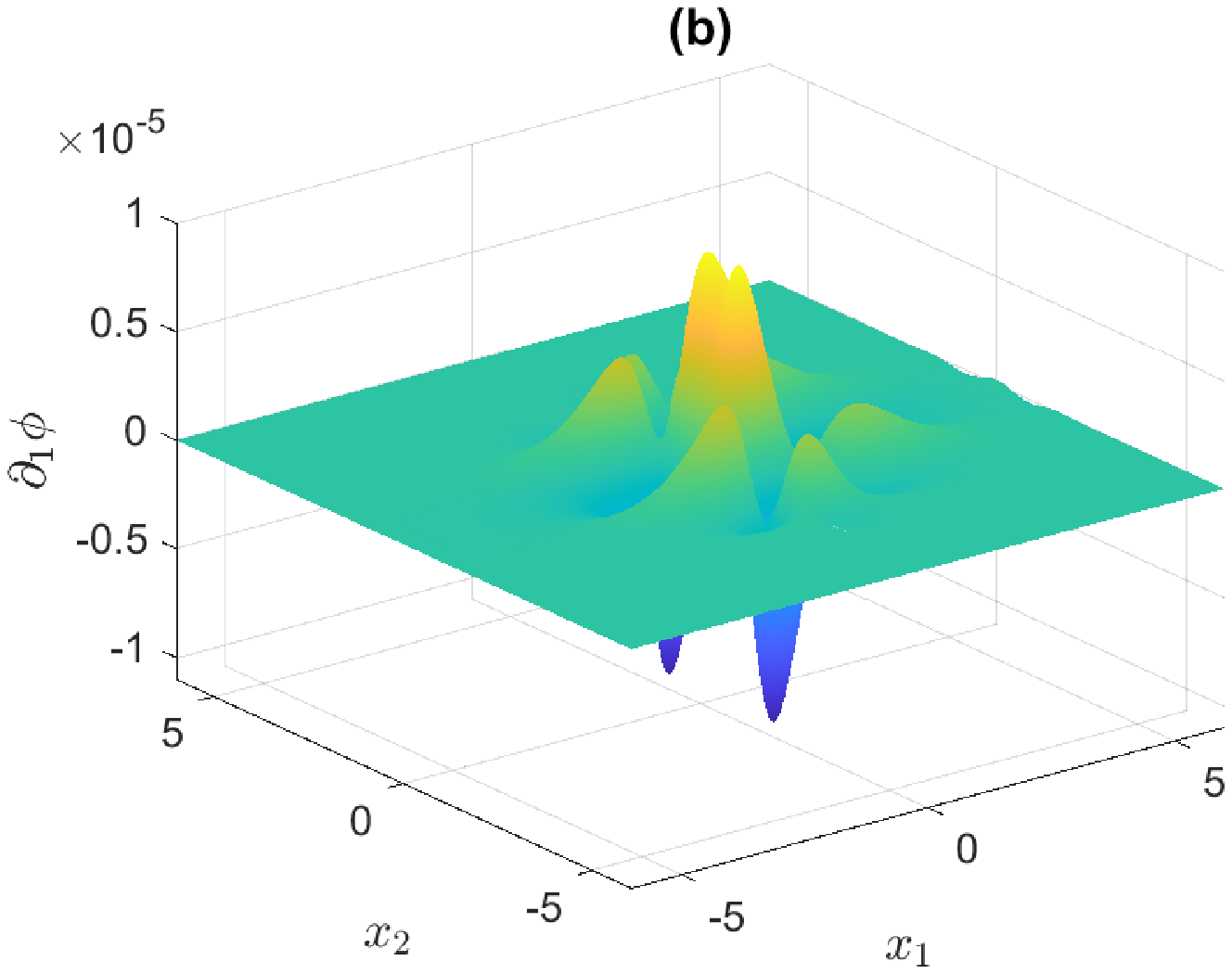}
	\end{subfigure}
	$\;$
	\begin{subfigure}[b]{0.38\linewidth}
		\includegraphics[width=\linewidth]{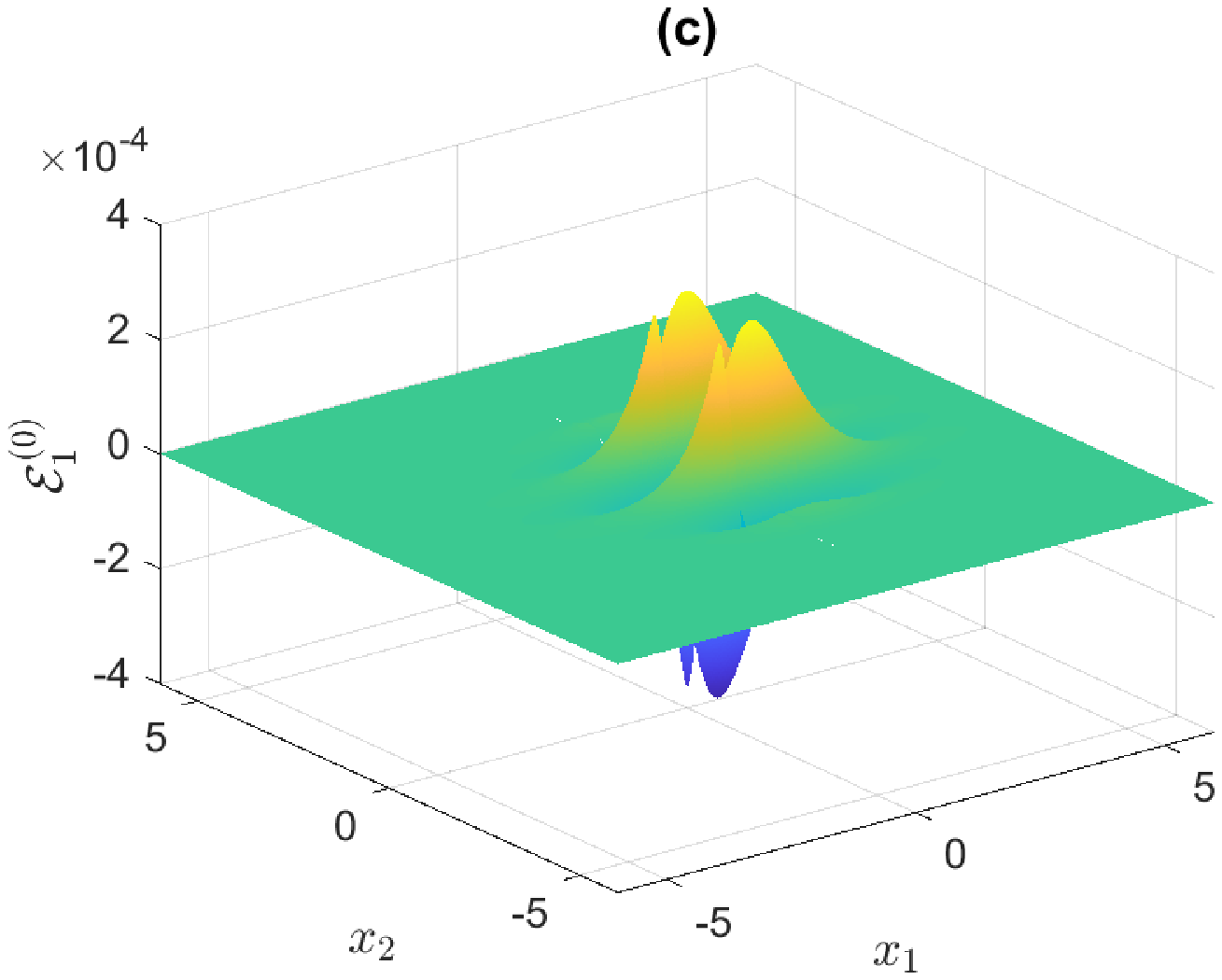}
	\end{subfigure}
	\caption{(a) The graph shows that numerically $\|\nabla \phi\|_2 = \mathcal{O}(\varepsilon^{3/2})$. (b) Plot of $\partial_1 \phi$ for $\varepsilon = 3 \cdot 10^{-4}$. (c) The first component of the corrected initial value $\mathcal{E}_1^{(0)} = U_{0,1} + \partial_1 \phi$ for $\varepsilon = 3 \cdot 10^{-4}$. Note the different scales in (b) and (c).}
	\label{Fig_fem}
\end{figure}

Finally, we study the $\varepsilon$-convergence of $\|\nabla \phi \|_2$. For the fixed step size $h = 0{.}005$ and $\varepsilon$ ranging from $10^{-4}$ to $10^{-3}$ we obtain the desired rate of convergence, see Figure \ref{Fig_fem} (a). For $\varepsilon = 3 \cdot 10^{-4}$ Figures \ref{Fig_fem} (b) and (c) show the first components the computed solutions $\nabla \phi$ and $\mathcal{E}^{(0)}$, respectively.

\bigskip
\section{Generalisations and open problems}\label{Section_openproblems}
Let us finally discuss possible generalisations and questions which are left unanswered by our analysis and which we believe to be of interest.

\medskip
\paragraph{Curved interfaces} The interface $\Gamma$ separating the two materials was chosen straight merely for reasons of simplicity. The results can be directly generalised also for a smooth curved interface separating $\R^n$ into two unbounded subsets.

\paragraph{Positivity of $\boldsymbol{\epsilon_1,\epsilon_f}$} The positivity imposed on the coefficients $\epsilon_1,\epsilon_f$ in (H0) is a common assumptions in the analysis literature because the operator involved in equation \eqref{eq} is then of elliptic type. We are not aware of results in the literature dealing with quasilinear operators with coefficients of varying sign.

However, in the electromagnetic context materials with negative material functions exist. For the nonlinear coefficient $\epsilon_f$ it is the case, e.g., for Galium compounds \cite{Sutherland}, and for metamaterials both $\epsilon_1$ and $\epsilon_f$ may be negative. We have also seen in Section \ref{Section_linear} that in the simplest case of a piecewise constant $\epsilon_1$, in order to obtain the existence of a TM eigenfunction, opposite signs of $\epsilon_1$ in the two subspaces need to be imposed, see the dispersion relation \eqref{disp_rel}. Hence, if the analytical results could be extended to the case of nonpositive coefficients, one may obtain a result similar to Proposition \ref{Prop_Anwendung} for a wider range of materials.

We believe that overcoming this sign problem would be of great interest both for the analysis and the physics application.
\medskip
\paragraph{Suboptimality of the asymptotics} As we already mentioned in Remark \ref{Remark_3/2}, we expect that the correction $\nabla\phi$ of $U_0$ is asymptotically of order $\cO(\varepsilon^{3/2})$ in the $H^1$-norm. For its $L^2$-norm, this behaviour has been also confirmed by our numerical test described in Section \ref{Section_numericaltest}, see Figure \ref{Fig_fem}(a). The suboptimality of the asymptotic estimate in Proposition \ref{Prop_Anwendung} ultimately relies on the term $\iii b$ in \eqref{estimate_phi}, since it produces the logarithmic term that we estimate by \eqref{b_log}. All other terms produced when applying Theorem \ref{Thm_phi} to our Maxwell context are of the expected order $\cO(\varepsilon^{3/2})$, see Section \ref{Section_ApplicationNL}. In order to get a sharp result we would need to get rid of $\iii b$ in \eqref{estimate_phi} and therefore to be able to deal with the product $\int b\phi$ differently than by \eqref{eq_L2_b}, i.e. than by means of Proposition \ref{OS_thm}.

Note that, once such an issue is solved for $\|\nabla\phi\|_2$, the problem of improving also the asymptotic estimates for the higher Sobolev norms does not arise. Indeed, the latter involve only terms which already bound $\|\nabla\phi\|_2$, and terms which are asymptotically of order $\cO(\varepsilon^{3/2})$.


\section*{Acknowledgements}
This research is supported by the German Research Foundation, DFG grant No. DO 1467/4-1.


\begin{thebibliography}{47}
	
	\bibitem{AG65} Aronszajn N., Gagliardo E. 
	{\em Interpolation spaces and interpolation methods.} Annali di Matematica 68 (1965), 51-117.
	
	\bibitem{crack} Atkinson C., Champion C.R.
	{\em A Mode III crack at the interface between two nonlinear materials.} Proceedings of the Royal Society of London. Series A, Mathematical and Physical Sciences 429, no. 1876 (1990), 247-257.
	
	\bibitem{AdAP} Azzollini A., d'Avenia P., Pomponio A.
	{\em Quasilinear elliptic equations in $\R^N$ via variational methods and Orlicz-Sobolev embeddings.} Calc. Var. Partial Differential Equations 49 (2014), no. 1-2, 197-213.
	
	\bibitem{BC} Badiale M., Citti G.
	{\em Concentration compactness principle and quasilinear elliptic equations in $\R^n$.} Comm. Partial Differential Equations 16 (1991), no. 11, 1795-1818.
	
	\bibitem{BF} Barile S., Figueiredo G.M.
	{\em Existence of least energy positive, negative and nodal solutions for a class of $p\&q$-problems with potentials vanishing at infinity.} J. Math. Anal. Appl. 427 (2015), no. 2, 1205-1233.
	
	\bibitem{BCS} Bartolo R., Candela A.M., Salvatore A.
	{\em On a class of superlinear $(p,q)$-Laplacian type equations on $\R^N$.} J. Math. Anal. Appl. 438 (2016), no. 1, 29-41.
	
	\bibitem{BBS} Bathory M., Bul\'{i}\v{c}ek M., Sou\v{c}ek O.
	{\em Existence and qualitative theory for nonlinear elliptic systems with a nonlinear interface condition used in electrochemistry.} Z. Angew. Math. Phys. 71 (2020), no. 3, Paper No. 74, 24 pp.

	\bibitem{Benci} Benci V., Fortunato D., Pisani L.
	{\em Soliton like solutions of a Lorentz invariant equation in dimension 3.} Rev. Math. Phys. 10 (1998), no. 3, 315-344. 
	
	\bibitem{BB} Benouhiba N., Belyacine Z.
	{\em A class of eigenvalue problems for the $(p,q)$-Laplacian in $\R^N$.} Int. J. Pure Appl. Math. 80 (2012), no. 5, 727-737.
	
	\bibitem{BHR} Berestycki H., Hamel F., Roques L.
	{\em Analysis of the periodically fragmented environment model. I. Species persistence.} J. Math. Biol. 51 (2005), no. 1, 75-113.

	\bibitem{Borsuk} Borsuk M.
	{\em Transmission problems for elliptic second-order equations in non-smooth domains.} Frontiers in Mathematics. Birkhäuser/Springer Basel AG, Basel, 2010.
		
	\bibitem{BDPW} Brown M., Dohnal T., Plum M., Wood I.
	{\em Spectrum of the Maxwell equations for the flat interface between homogeneous dispersive media.} In preparation.
	
	\bibitem{CZ} Chen Zh., Zou J.
	{\em Finite element methods and their convergence for elliptic and parabolic interface problems.} Numer. Math. 79 (1998), no. 2, 175-202.
	
	\bibitem{bio} Colli Franzone P., Guerri L., Rovida, S.
	{\em Wavefront propagation in an activation model of the anisotropic cardiac tissue: asymptotic analysis and numerical simulations.} J. Math. Biol. 28 (1990), 121-176.
	
	\bibitem{CDN} Costabel M., Dauge M., Nicaise S.
	{\em Singularities of Maxwell interface problems.} M2AN Math. Model. Numer. Anal. 33 (1999), no. 3, 627-649.
	
	\bibitem{DR} Dohnal T., Rudolf D.
	{\em NLS approximation for wavepackets in periodic cubically nonlinear wave problems in $\R^d$.} Appl. Anal. 99 (2020), no. 10, 1685-1723.
	
	\bibitem{EFK} Ebmeyer C., Frehse J., Kassmann M.
	{\em Boundary regularity for nonlinear elliptic systems: applications to the transmission problem.} Geometric Analysis and Nonlinear Partial Differential Equations (S.Hildebrandt and H.Karcher, eds.). Springer Verlag 2002, pp. 505-517. 
	
	\bibitem{Feynman} Feynman R.P., Leighton R.B., Sands M.
	{\em The Feynman lectures on physics. Vol. 2: Mainly electromagnetism and matter.} Addison-Wesley Publishing Co., Inc., Reading, Mass.-London 1964
	
	\bibitem{Galdi} Galdi G.P.
	{\em An introduction to the mathematical theory of the Navier-Stokes equations. Steady-state problems.} Second edition. Springer Monographs in Mathematics. Springer, New York, 2011.
	
	\bibitem{Giusti} Giusti E.
	{\em Direct methods in the calculus of variations.} World Scientific Publishing Co., Inc., River Edge, NJ, 2003.
	
	\bibitem{KSM} Kirrmann P., Schneider G., Mielke A.
	{\em The validity of modulation equations for extended systems with cubic nonlinearities.} Proc. Roy. Soc. Edinburgh Sect. A 122 (1992), no. 1-2, 85-91.
	
	\bibitem{Knees} Knees D.
	{\em On the regularity of weak solutions of quasi-linear elliptic transmission problems on polyhedral domains.} Z. Anal. Anwendungen 23 (2004), no. 3, 509-546.
		
	\bibitem{KS} Knees D., S\"{a}ndig A.-M.
	{\em Regularity of elastic fields in composites.} Multifield problems in solid and fluid mechanics, 331-360, Lect. Notes Appl. Comput. Mech., 28, Springer, Berlin, 2006.
	
	\bibitem{Kuro} Kurokawa T.
	{\em Riesz potentials, higher Riesz transforms and Beppo Levi spaces.} Hiroshima Math. J. 18 (1988), no. 3, 541-597.
	
	\bibitem{KL} Kutev N., Lions P.-L.
	{\em Nonlinear second-order elliptic equations with jump discontinuous coefficients. Part I, quasilinear equations.} Differential and Integral Equations 5 (1992), no. 6, 1201-1217.
	
	\bibitem{LSU} Ladyzhenskaja O.A., Solonnikov V.A., Ural'tseva N.N.
	{\em Linear and quasilinear equations of parabolic type.} (Russian)
	Translated from the Russian by S. Smith Translations of Mathematical Monographs, Vol. 23 American Mathematical Society, Providence, R.I. 1968.
	
	\bibitem{LRU} Ladyzhenskaja O.A., Rivkind V.Ja., Ural'tseva N.N.
	{\em Solvability of diffraction problems in the classical sense.} (Russian) Trudy Mat. Inst. Steklov. 92 (1966), 116-146.
	
	\bibitem{LU} Ladyzhenskaya O.A., Ural'tseva N.N.
	{\em Linear and quasilinear elliptic equations.} Academic Press, New York-London 1968.
	
	\bibitem{LS} Li G., Seshadri S.R. 
	{\em Weakly nonlinear surface polariton.} J. Opt. Soc. Am. B 6 (1989) 1125-1137.
	
	\bibitem{LL} Li G., Liang X.
	{\em The existence of nontrivial solutions to nonlinear elliptic equation of $(p,q)$-Laplacian type on $\R^N$.} Nonlinear Anal. 71 (2009), no. 5-6, 2316-2334.

	\bibitem{Liu} Liu W.B.
	{\em Degenerate quasilinear elliptic equations arising from bimaterial problems in elastic-plastic mechanics.} Nonlinear Anal. 35 (1999), no. 4, Ser. A: Theory Methods, 517-529.
	
	\bibitem{Lorenzi} Lorenzi A.
	{\em On elliptic equations with piecewise constant coefficients II} Annali della Scuola Normale Superiore di Pisa, Classe di Scienze 3e s\'{e}rie, tome 26, no. 4 (1972), 839-870.
	
	\bibitem{Mar} Marcellini P.
	{\em Regularity and existence of solutions of elliptic equations with $p,q$-growth conditions.} J. Differential Equations 90 (1991), no. 1, 1-30.
	
	\bibitem{Mar2} Marcellini P.
	{\em Growth conditions and regularity for weak solutions to nonlinear elliptic pdes.} J. Math. Anal. Appl. (2020). https://doi.org/10.1016/j.jmaa.2020.124408.
	
	\bibitem{MS} Marini A., Skryabin D.V.
	{\em Ginzburg-Landau equation bound to the metal-dielectric interface and transverse nonlinear optics with amplified plasmon polaritons.} Phys. Rev. A 81 033850 (2010).
	
	\bibitem{M} Mingione G.
	{\em Regularity of minima: an invitation to the dark side of the calculus of variations.} Appl. Math. 51 (2006), no. 4, 355-426.
	
	\bibitem{MR} Mingione G., R\u{a}dulescu V.
	{\em Recent developments in problems with nonstandard growth and nonuniform ellipticity.} J. Math. Anal. Appl. (2021). https://doi.org/10.1016/j.jmaa.2021.125197.
	
	\bibitem{MR_thermo} Mu\~{n}oz Rivera J.E., Portillo Oquendo H.
	{\em A transmission problem for thermoelastic plates.} Quart. Appl. Math. 62 (2004), no. 2, 273-293.
	
	\bibitem{NS} Nicaise S., Sändig A.-M.
	{\em Transmission problems for the Laplace and elasticity operators.} Regularity and boundary integral formulation. Math. Meth. Appl. Sci. 9 (1999), 855-898.
	
	\bibitem{Oleinik} Ole\u{\i}nik O.A.
	{\em Boundary-value problems for linear equations of elliptic parabolic type with discontinuous coefficients.} (Russian) Izv. Akad. Nauk SSSR Ser. Mat. 25 (1961), 3-20.
	
	\bibitem{OS} Ortner C., S\"{u}li E. {\em A note on linear elliptic systems on $\R^d$}. Arxiv preprint (2012), arXiv:1202.3970.
	
	\bibitem{S} Simader Ch.G. {\em The weak Dirichlet and Neumann problem for the Laplacian in $L^q$ for bounded and exterior domains. Applications.}  In: Krbec M., Kufner A., Opic B., Rákosník J. (eds) Nonlinear Analysis, Function Spaces and Applications Vol. 4. TEUBNER-TEXTE zur Mathematik. Vieweg+Teubner Verlag, Wiesbaden.
	
	\bibitem{SS} Schnaubelt R., Spitz M.
	{\em Local wellposedness of quasilinear Maxwell equations with conservative interface conditions.} Arxiv preprint (2018), arXiv:1811.08714.
	
    \bibitem{Sutherland} Sutherland R.L. 
    {\em Handbook of Nonlinear Optics.} Optical Science and Engineering, Taylor \& Francis 2003.
		
	\bibitem{TL2} Tan Q.J., Leng Z.J. {\em The method of upper and lower solutions for diffraction problems of quasilinear elliptic reaction-diffusion systems.} J. Math. Anal. Appl. 380 (2011), no. 1, 363-376.
	
	\bibitem{Daniel} Tietz D.P.
	{\em Justification of the Nonlinear Schr\"odinger Equation for Interface Wave Packets in Maxwell's Equations with 2D Localization.} Ph.D. Dissertation, Martin-Luther Universit\"at Halle-Wittenberg. In preparation.
	
	\bibitem{TT} Truesdell C., Toupin R. {\em The classical field theories.} With an appendix on tensor fields by J. L. Ericksen. In: Fl\"{u}gge, S. (ed.) Handbuch der Physik, Bd. III/1, pp. 226-793, appendix, pp. 794-858. Springer, Berlin (1960).
	
	\bibitem{Wang} Wang, W-Ch.
	{\em A jump condition capturing finite difference scheme for elliptic interface problems.} SIAM J. Sci. Comput. 25 (2004), no. 5, 1479-1496.
	
	\bibitem{Whitney} Whitney H.
	{\em Analytic extensions of differentiable functions defined in closed sets.} Trans. Amer. Math. Soc. 36 (1934), 63-89.
	
\end{thebibliography}
\end{document}